\documentclass[reqno]{amsart}

\usepackage{amsmath}
\usepackage{amssymb}
\usepackage{amsthm}
\usepackage{mathtools}
\usepackage{enumitem}
\setenumerate{label={\normalfont (\roman*)}}

\usepackage[utf8]{inputenc}
\usepackage[T1]{fontenc}
\usepackage{lmodern}
\usepackage[babel]{microtype}
\usepackage[english]{babel}
\usepackage{relsize}

\usepackage{graphicx}
\usepackage{subcaption}

\usepackage{geometry}
\usepackage{xcolor} 	
\usepackage{hyperref}
\hypersetup{
	colorlinks,
	linkcolor={red!60!black},
	citecolor={green!60!black},
	urlcolor={blue!60!black},
}
\usepackage[abbrev, msc-links]{amsrefs}
\usepackage[nameinlink, capitalise, noabbrev]{cleveref}
\crefformat{enumi}{#2#1#3}
\crefformat{equation}{#2(#1)#3}
\usepackage{doi}

\renewcommand{\PrintDOI}[1]{\doi{#1}}

\let\setminus=\smallsetminus

\renewcommand{\subset}{\subseteq}
\renewcommand{\supset}{\supseteq}
\renewcommand{\leq}{\leqslant}
\renewcommand{\geq}{\geqslant}
\renewcommand{\ge}{\geq}
\renewcommand{\le}{\leq}

\newtheorem{theorem}{Theorem}[section]

\newtheorem{lemma}[theorem]{Lemma}
\newtheorem{keylemma}[theorem]{Key Lemma}

\newtheorem{observation}[theorem]{Observation}

\newtheorem{mainresult}{Theorem} 

\theoremstyle{definition}
\newtheorem{example}[theorem]{Example}

\theoremstyle{remark}

\newtheorem*{ack}{Acknowledgement}

\let\theta=\vartheta
\let\rho=\varrho
\let\phi=\varphi
\def\NN{\mathbb N}

\def\QQ{\mathbb Q}

\def\cG{{\mathcal G}}
\def\cH{{\mathcal H}}

\def\cL{{\mathcal L}}
\def\cM{{\mathcal M}}

\def\cP{{\mathcal P}}
\def\cQ{{\mathcal Q}}

\def\cU{{\mathcal U}}

\def\cX{{\mathcal X}}

\newcommand{\set}[1]{{\lbrace {#1} \rbrace}}
\newcommand{\oset}[1]{( {#1} )}
\newcommand{\bigset}[1]{{\left\lbrace {#1} \right\rbrace}}

\usepackage{stmaryrd}
\usepackage{letterswitharrows}

\renewcommand{\leq}{\leqslant}
\renewcommand{\geq}{\geqslant}

\newcommand{\Kaleph}{K^{\aleph_0}}
\newcommand{\tcd}{tree-cut de\-com\-po\-si\-tion}
\newcommand{\FG}{F}
\newcommand{\hFG}{\breve{F}}

\newcommand{\Pigraph}{$\Pi$-graph}

\newcommand{\free}{free}
\newcommand{\ccon}{compound-connected}
\newcommand{\csep}{compound-separation}
\newcommand{\cseparator}{separator}

\newcommand{\unit}{unitary}

\newcommand{\ccv}{compound-cutvertex}
\newcommand{\ccvs}{compound-cutvertices}
\newcommand{\partordrest}[2]{#1\!\upharpoonright\!#2}

\begin{document}
\setlength{\fboxsep}{0pt}
\setlength{\fboxrule}{.1pt}

\title{The immersion-minimal infinitely edge-connected graph}

\author{Paul Knappe}
\address{Universität Hamburg, Department of Mathematics, Bundesstraße 55 (Geomatikum), 20146 Hamburg, Germany}
\email{paul.knappe@uni-hamburg.de}
\author{Jan Kurkofka}
\address{University of Birmingham, Birmingham, UK}
\email{j.kurkofka@bham.ac.uk}


\begin{abstract}
We show that there is a unique immersion-minimal infinitely edge-connected graph: every such graph contains the halved Farey graph, which is itself infinitely edge-connected, as an immersion minor.

By contrast, any minimal list of infinitely edge-connected graphs represented in all such graphs as topological minors must be uncountable.
\end{abstract}

\vspace*{-3cm}
\maketitle

\vspace*{-.7cm}

\section{Introduction}

\begin{figure}[h]
    \centering
    \includegraphics[width=.35\textwidth]{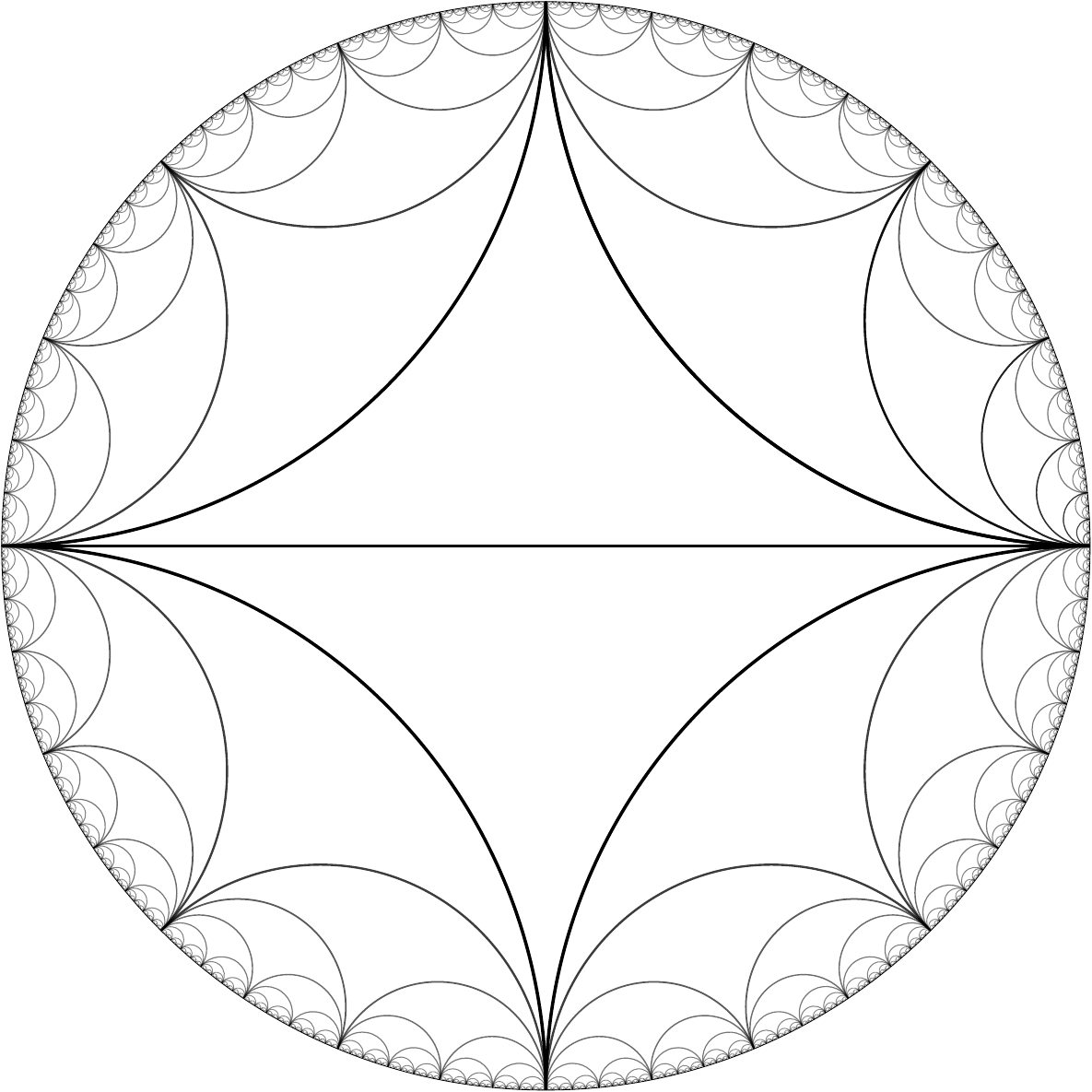}
    \caption{The Farey graph}
    \label{fig:FareyGraph}

        
\end{figure}

\noindent The Farey graph, shown in \cref{fig:FareyGraph} and surveyed in \cites{OfficeHoursGroupTheory,hatcher2017topology}, plays a role in a number of mathematical fields ranging from group theory and number theory to geometry and dynamics~\cite{OfficeHoursGroupTheory}.
Curiously, graph theory has not been among these until very recently, when it was shown in~\cite{kurkofka2020infinitelyedgecongraphcontaingfareygraphorT} that the Farey graph plays a central role in graph theory too: 

\begin{theorem}\label{thm:originalFarey}
The Farey graph is one of two infinitely edge-connected graphs such that every infinitely edge-connected graph contains at least one of the two as a minor.
\end{theorem}

Naturally, this result raises the question of whether there exist similar infinitely edge-connected graphs for graph relations other than the minor relation which determine infinite edge-connectivity by forming minimal lists in this way.
In this paper, we address two graph relations that harmonise particularly well with edge-connectivity: the topological minor relation and the immersion relation.

A \emph{weak immersion} of a graph~$H$ in a graph~$G$ is a map~$\alpha$ with domain~$V(H)\sqcup E(H)$ that embeds~$V(H)$ into~$V(G)$ and sends every edge~$uv\in H$ to an~$\alpha(u)$--$\alpha(v)$ path in~$G$ which is edge-disjoint from every other such path.
The map~$\alpha$ is a \emph{strong immersion} of~$H$ in~$G$ if additionally all paths~$\alpha(e)$ for~$e\in E(H)$ have no internal vertices in~$\alpha[V(H)]$.
The vertices of~$G$ that lie in the image~$\alpha[V(H)]$ are the \emph{branch vertices} of this immersion.

We say that~$H$ is \emph{strongly immersed} in~$G$, or that~$H$ is a \emph{strong immersion minor} of~$G$, if there is a strong immersion of~$H$ in~$G$.
Similarly, we define \emph{weakly immersed} and \emph{weak immersion minor}.
Robertson and Seymour showed that the weak immersion relation well-quasi-orders the finite graphs, just like the minor relation, and they believe that so does the strong immersion relation~\cite{GMXIII}*{§1}.
In this paper, we will focus on strong immersions. 
For brevity, we will often refer to strong immersions simply as \emph{immersions}.

Any infinitely edge-connected graphs that form a minimal list as discussed earlier must be countable, because in every infinitely edge-connected graph we can greedily find a countable infinitely edge-connected subgraph.
The countable graphs, however, are not known to be well-quasi-ordered by the minor relation or either of the immersion relations.
It is therefore not clear that any immersion-minimal set of infinitely edge-connected graphs must be finite, nor even that such a minimal set exists.

\begin{figure}[t]
\centering
\begin{minipage}{.49\textwidth}
  \centering
  \includegraphics[width=.7\linewidth]{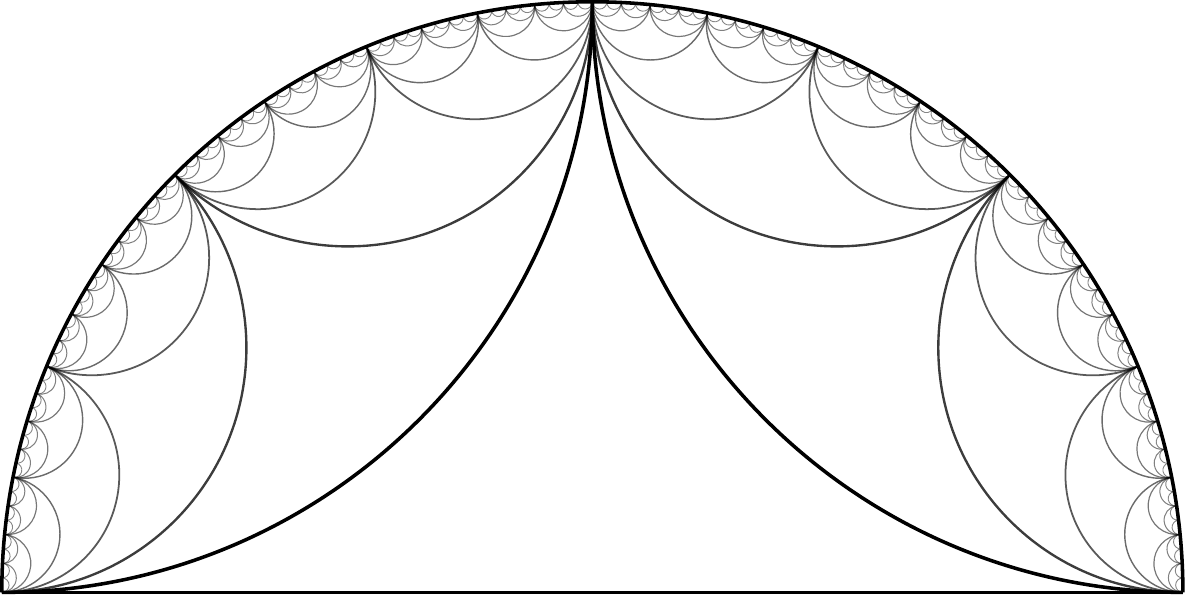}
  \captionof{figure}{The halved Farey graph}
  \label{fig:halvedFareyGraph}
\end{minipage}%
\begin{minipage}{.49\textwidth}
  \centering
  \captionsetup{width=.9\linewidth}
  \includegraphics[width=.7\linewidth]{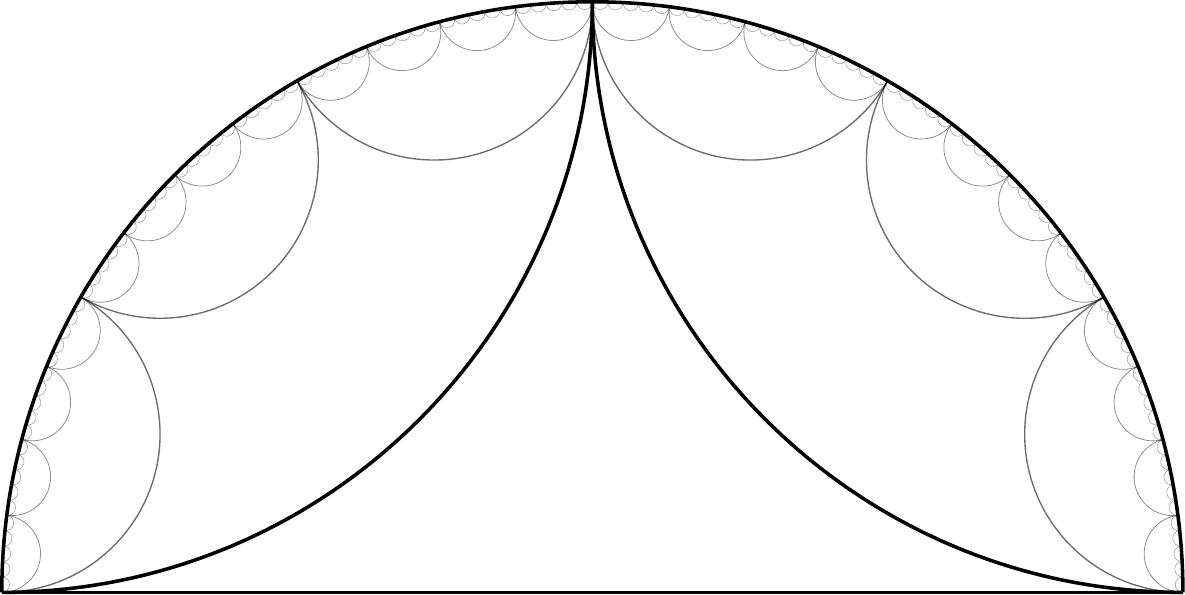}
  \captionof{figure}{A generalised halved Farey graph}
  \label{fig:genfareygraph}
\end{minipage}

    
\end{figure}

A greedy argument shows that every infinitely edge-connected graph contains the countably infinite complete graph by weak immersion.
So for weak immersion we have a minimal list formed by this graph alone.
But this is no longer true if we replace `weak immersion' with `immersion'.
Indeed, we show that the countably infinite complete graph is not immersed in the Farey graph (\cref{thm:usualcandidatesandimmersion}~\cref{item:infinitecomplete}).
Then, does every infinitely edge-connected graph contain the Farey graph by immersion?
Perhaps surprisingly, the answer is no: the Farey graph is not an immersion minor of the \emph{halved} Farey graph shown in \cref{fig:halvedFareyGraph} (\cref{thm:usualcandidatesandimmersion}~\cref{item:Farey}), although the halved Farey graph is infinitely edge-connected. 
As~our main result we show that the halved Farey graph is immersed in every infinitely edge-connected graph, and hence forms the desired list for immersion, again all by itself.
Two graphs that are immersed in each other are called \emph{immersion-equivalent}.

\begin{mainresult} \label{mainthm:halvedFareygraphistypical}
Up to immersion equivalence, the halved Farey graph is the unique infinitely edge-connected graph that is immersed in every infinitely edge-connected graph.
\end{mainresult}

Neither the minor relation nor the immersion relation is stronger than the other.
Therefore, \cref{thm:originalFarey} and \cref{mainthm:halvedFareygraphistypical} do not obviously imply each other.
The topological minor relation, however, is stronger than both the minor relation and the immersion relation.
This raises the question of whether the two theorems can be unified by finding a short list of topological-minor-minimal infinitely edge-connected graphs.
Our second result, \cref{mainthm:anytypicalsetisbig} below, shows that this is impossible.
Hence, \cref{thm:originalFarey} and~\cref{mainthm:halvedFareygraphistypical} are best possible in this sense.

Let~$\cG$ be a class of graphs and~$\le$ a relation on~$\cG$.
We say that a class~$\cH\subset\cG$ is \emph{typical} for~$\cG$ with regard to~$\le$ if for every graph~$G \in \cG$ there exists a graph~$H \in \cH$ with~$H \le G$.
The relation of `being typical for' with regard to~$\le$ is transitive on the subclasses of~$\cG$.
We recall that in every infinitely edge-connected graph we can greedily find a countable infinitely edge-connected subgraph.
Combining this with the fact that every countable graph is isomorphic to a subgraph of~$K^{\aleph_0}$, and using the transitivity of the `being typical for'-relation, we find that the class of infinitely edge-connected graphs includes typical sets of graphs with regard to the topological minor relation that are no larger than the continuum.
However, we show that none of these typical sets is countable, let alone finite.
Recall that a graph is \emph{outerplanar} if it has a drawing in which every vertex lies on the unit circle
and every edge is contained in the unit disc. Indeed:

\begin{mainresult} \label{mainthm:anytypicalsetisbig}
Every set of graphs that is typical for the infinitely edge-connected graphs with regard to the topological minor relation, or even just for the outerplanar ones, is uncountable.
\end{mainresult}

\noindent Our proof of \cref{mainthm:anytypicalsetisbig} builds on a construction of generalised halved Farey graphs. An instance of a generalised halved Farey graph is shown in \cref{fig:genfareygraph}.

This paper is organised as follows.
In \cref{sec:prelims} we introduce the tools and terminology that we need.
In \cref{sec:genfareygraph} we define the generalised halved Farey graphs and prove \cref{mainthm:anytypicalsetisbig}.
In \cref{sec:StrongImmersionProof} we prove \cref{mainthm:halvedFareygraphistypical}.

\section{Tools and terminology}\label{sec:prelims}

\noindent We use the notation of Diestel's book~\cite{DiestelBook}. 
Recall that a~non-trivial path~$P$ is an~\emph{$A$-path} for a set~$A$ of vertices if~$P$ has its endvertices but no inner vertex in~$A$.
Given a graph~$H$, we call~$P$ an \emph{$H$-path} if~$P$ is non-trivial and meets~$H$ exactly in its endvertices.
An $H$-path is a $V(H)$-path; the converse is true unless the $V(H)$-path has just one edge and this edge is in~$H$.
We write~$G[X]$ for the subgraph of~$G$ induced by the vertex set~$X$.
Given a path~$P$ that contains two vertices~$u$ and~$v$, we write~$uPv$ for the subpath of~$P$ from~$u$ to~$v$.

Whenever an~$x$--$y$ path~$P$ is introduced, we denote by~$\leq_{P}$ the linear order on its vertices given by the way that~$P$ directed from~$x$ to~$y$ traverses them.
For a partial order~$\cL=(L, \leq_L)$ and a subset~$M$ of~$L$, we write~$\partordrest{\cL}{M}$ for the partial order obtained from~$\cL$ by restricting it to~$M$.

\subsection{Farey graph}
The \emph{Farey graph}~$F$ is the graph on~$\QQ\cup\{\infty\}$ in which two rational numbers~$a/b$ and~$c/d$ in lowest terms (allowing also~$\infty=(\pm 1)/0$) form an edge if and only if~$\det\bigl( \begin{smallmatrix}a & c\\ b & d\end{smallmatrix}\bigr)=\pm 1$, cf.~\cite{OfficeHoursGroupTheory}.
In this paper we do not distinguish between the Farey graph and the graphs that are isomorphic to it.
For our graph-theoretic proofs it will be more convenient to work with the following purely combinatorial definition of the Farey graph that is indicated in~\cite{OfficeHoursGroupTheory} and~\cite{hatcher2017topology}.

The \emph{halved Farey graph~$\hFG_0$ of order~$0$} is a~$K^2$ with its sole edge coloured blue.
Inductively, the \emph{halved Farey graph~$\hFG_{n+1}$ of order~$n+1$} is the edge-coloured graph that is obtained from~$\hFG_n$ by adding a new vertex~$v_e$ for every blue edge~$e \in \hFG_n$, joining every~$v_e$ precisely to the endvertices of~$e$ by two blue edges, and recolouring all the edges of~$\hFG_{n+1}$ belonging to~$\hFG_n$ black.
The \emph{halved Farey graph}~$\hFG:=\bigcup_{n\in\NN}\hFG_n$ is the union of all these~$\hFG_n$ without their edge-colourings (cf. \cref{fig:halvedFareyGraph}), and the \emph{Farey graph} is the union~$F=G_1\cup G_2$ of two copies~$G_1,G_2$ of the halved Farey graph such that~$G_1\cap G_2=\hFG_0$ (cf. \cref{fig:FareyGraph}).

The (halved) Farey graph and infinite edge-connectivity have been studied and used in~\cites{kurkofka2020infinitelyedgecongraphcontaingfareygraphorT,kurkofka2020fareydeterminedbyconnectivity,kurkofka2019ubiquityandfarey,Grid}.

\subsection{Grain lines}

Suppose that~$x$ and~$y$ are two vertices in a graph~$G$ such that no finite set of edges separates~$x$ and~$y$ in~$G$.
Then we greedily find a sequence~$\oset{P_n : n \in \NN}$ of infinitely many pairwise edge-disjoint~$x$--$y$ paths in~$G$.
Since these paths are only edge-disjoint, they can meet in vertices other than~$x$ and~$y$.
Let us say that two~$x$--$y$ paths are \emph{order-compatible} if they traverse their common vertices in the same order.
Is it always possible to choose the paths~$P_n$ so that they are pairwise order-compatible?
Perhaps surprisingly, the answer is no:
in~\cite{kurkofka2019ubiquityandfarey}, a countable planar graph is constructed that is infinitely edge-connected, but which does not contain infinitely many edge-disjoint pairwise order-compatible paths between any two of its vertices.

Fortunately, not all is lost.
While we cannot always choose the paths~$P_n$ so that they are pairwise order-compatible, we can always choose them so that they satisfy a slightly weaker form of order-compatibility which is still strong enough for our purpose.
Roughly speaking, we will be able to choose the paths~$P_n$ so that they induce a linear order on their limit.
This limit will be the set of all vertices that eventually appear on all paths~$P_n$.
If two vertices~$u$ and~$v$ are in the limit, there will be a first path~$P_n$ which contains both~$u$ and~$v$.
The~$x$--$y$ path~$P_n$ linearly orders its vertex set from~$x$ to~$y$; in particular, it orders~$u$ and~$v$.
This order on~$u$ and~$v$ might disagree with a later path that is not order-compatible with~$P_n$.
However, we shall achieve that all later paths will order~$u$ and~$v$ in the same way, and we will use this ordering of~$u$ and~$v$ in our limit instead of the ordering induced by~$P_n$.
This informal idea has been formalised as `grain lines' in~\cite{kurkofka2020fareydeterminedbyconnectivity}, whose definition we recall now.

An \emph{$x$--$y$ grain line} between two distinct vertices~$x$ and~$y$ is an ordered pair~$(\cL, \cP)$ where~$\cL = (L, \leq_L)$ is a linear order with~$\min_{\leq_L} L = x$ and~$\max_{\leq_L} L = y$, and~$\cP = \oset{P_n : n \in \NN}$ is a sequence of pairwise edge-disjoint~$x$--$y$ paths~$P_n$ such that the following conditions are satisfied:
\begin{enumerate}[label=(GL\arabic*)]
    \item \label{prop:Lislimit} $L = \set{v : \set{n \in \NN : v \in V(P_n)} \text{ is a final segment of } \NN}$;
    \item \label{prop:externallydisjoint} if a vertex of a path~$P_n$ is not contained in~$L$, then it is not a vertex of any other path~$P_m$;
    \item \label{prop:samelinearorder} for every~$n \geq 1$, the linear order~$\leq_{P_n}$ given by~$P_n$ and~$\leq_L$ induce the same linear order on the vertex set~$L_{<n}$, where we set~$L_{<n} := L \cap \bigcup_{m < n} V(P_m)$.
\end{enumerate}

\begin{example} \label{example:halvedfareyisgrain}
The halved Farey graph defines a grain line, as follows.
Let~$x$ and~$y$ be the two vertices of~$\hFG_0$.
For every~$n\in\NN$, let~$P_n$ be the blue Hamilton path of~$\hFG_n$, and let us view each~$P_n$ as an~$x$--$y$ path in~$\hFG$.
Then letting~$L:=V(\hFG)$, $\le_L:=\bigcup_{n\in\NN}\le_{P_n}$ and~$\cP:=(P_n: n\in\NN)$ results in an~$x$--$y$ grain line.
\end{example}

\begin{lemma} \label{thm:existenceofgrainlines}
    Let~$x$ and~$y$ be any two distinct vertices of a graph~$G$, and let~$\cQ$ be any set of infinitely many pairwise edge-disjoint~$x$--$y$ paths in~$G$.
    Then there exists an~$x$--$y$ grain line~$(\cL,\cP)$ in~$G$
    such that all the paths in~$\cP$ are in~$\cQ$.
\end{lemma}

\begin{proof}
    The proof of \cite{kurkofka2020fareydeterminedbyconnectivity}*{Theorem 5.4} shows this.
\end{proof}

Whenever a grain line is introduced as~$(\cL,\cP)$, we tacitly assume that~$\cL = (L, \leq_L)$ and~$\cP = \oset{P_n : n \in \NN}$.
We write~$\bigcup \cP := \bigcup_{n \in \NN} P_n$ for the graph \emph{defined} by the grain line~$(\cL, \cP)$.
A \emph{$\cP$-segment} is a subpath~$uP_dv$ of some path~$P_d$ in~$\cP$ with~$d \geq 1$, if~$u$ and~$v$ are in~$L_{<d}$ and~$v$ is the successor of~$u$ in~$\partordrest{\cL}{L_{<d}}$.
We follow the convention that~$P_0$ also is a $\cP$-segment.
We refer to~$d$ as the \emph{$\cP$-depth} of the~$\cP$-segment~$uP_dv$, and the $\cP$-depth of~$P_0$ is~0.

\begin{example}
    Every~$\cP$-segment of a grain line has at least one edge. It is possible for~$\cP$-segments to have only one edge: for example, if all paths~$P_n\in\cP$ are internally disjoint, and all paths have two edges except one path which has exactly one edge.
\end{example}

We introduce the concept of depth in a sequence~$\cP$ of paths.
The path~$P_d$ is the path in \emph{$\cP$-depth}~$d$.
For every~$d \in \NN$, we abbreviate the sequence~$\oset{P_n : n \geq d}$ of paths in depth at least~$d$ as~$\cP_{\geq d}$.
Similarly, we define~$\cP_{> d}$, $\cP_{\leq d}$ and~$\cP_{<d}$.
The \emph{$\cP$-depth} of a vertex~$v\in \bigcup \cP$ is defined as~$\min\set{n : v \in V(P_n)}$.
The \emph{$\cP$-depth} of an edge~$e \in \bigcup \cP$ is~$\min\set{n : e \in E(P_n)}$.
The \emph{depth} of a vertex or an edge of~$\bigcup \cP$ in~$\cP$ is its~$\cP$-depth.

Let~$(\cL,\cP)$ be a grain line.
We remark that~$L_{<d}$ is the set of vertices in~$L$ whose depth is less than~$d$ in~$\cP$, and the depth of an edge in~$\cP$ is at least the depth of its endvertices in~$\cP$.
For two vertices~$u, v \in L$ with~$u <_L v$, we write~$u\cP v$ for the subsequence~$\oset{uP_nv : n > d}$ where~$d$ is the maximum of the depths of~$u$ and~$v$ in~$\cP$.

The following structural properties have been introduced in \cite{kurkofka2020fareydeterminedbyconnectivity}.
A grain line~$(\cL, \cP)$ is \emph{wild} if~$\cL$ is order-isomorphic to~$\QQ \cap [0,1]$.
It is \emph{wildly presented} if, for every~$n \geq 1$, whenever~$u <_L v$ are elements of~$L_{<n}$ then an internal vertex of~$uP_nv$ is in the interval~$(u,v)_\cL$.

\section{Typical sets with regard to the topological minor relation}\label{sec:genfareygraph}

\noindent In this section, we show the following generalisation of \cref{mainthm:anytypicalsetisbig}. Missing definitions follow.

\begin{theorem} \label{thm:countablesetsnottypical}
    For every countable set of infinitely edge-connected graphs, there exists an outerplanar infinitely edge-connected graph that contains none of them as a topological minor.
    Moreover, there is such a graph which is a generalised halved Farey graph.
\end{theorem}

A \emph{generalised halved Farey graph of order~$0$} is a non-trivial path with its edges coloured blue.
Inductively, a \emph{generalised halved Farey graph of order~$n+1$} is an edge-coloured graph that is obtained from a generalised Farey graph~$G_n$ of order~$n$ by adding, for every blue edge~$e=uv \in G_n$, a blue-coloured~$u$--$v$ path~$P_e$ of length at least two, which is internally disjoint from~$G_n$ and every other~$P_{e'}$, and recolouring all the edges of~$G_n$ black in~$G_{n+1}$.
Let~$\oset{G_n : n \in \NN}$ be any sequence obtained by this construction.
Then the union of all these~$G_n$ without their edge-colourings is a \emph{generalised halved Farey graph}.
We remark that every generalised halved Farey graph is outerplanar.

Just like the halved Farey graph in \cref{example:halvedfareyisgrain}, we may interpret a generalised halved Farey graph as a grain line.
This grain line then satisfies the following stronger versions of \cref{prop:externallydisjoint} and \cref{prop:samelinearorder}:
\begin{enumerate}[label=(GL\arabic*')] \setcounter{enumi}{1}
    \item \label{GL2'} $L=\bigcup_{n\in\NN}V(P_n)$;
    \item \label{GL3'} ${\le_L}=\bigcup_{n\in\NN}{\le_{P_n}}$.
\end{enumerate}
It is immediate to see that the converse is also true:

\begin{lemma}\label{lem:genhalvedfareygraphisgrainline}
    A graph~$G$ is a generalised halved Farey graph if and only if~$G$ is defined by a grain line~$(\cL, \cP)$ which satisfies \cref{GL2'} and \cref{GL3'}.
    \qed
\end{lemma}

The generalised halved Farey graphs which we will construct in the proof of \cref{thm:countablesetsnottypical} are of the following type.
Let~$\ell \colon \NN \to \NN$ be a function such that~$\ell(0) \geq 1$ and, for every~$n\geq 1$, $\ell(n) \geq 2$.
Then~$\ell$ \emph{induces} an (up to isomorphism) unique grain line~$(\cL,\cP)(\ell)$ which satisfies \cref{GL2'}, \cref{GL3'} and, for every~$n \in \NN$, every~$\cP$-segment in~$\cP$-depth~$n$ is a path of length~$\ell(n)$.
The \emph{generalised halved Farey graph~$\hFG\oset{\ell}$ induced by~$\ell$} is the graph defined by~$(\cL,\cP)(\ell)$.
Note that by \cref{lem:genhalvedfareygraphisgrainline}, $\hFG\oset{\ell}$ is a generalised halved Farey graph.

\begin{example} 
The halved Farey graph is the generalised halved Farey graph induced by~$\ell$ with~$\ell(0) := 1$ and~$\ell(n) := 2$ for all natural numbers~$n \geq 1$.
The generalised halved Farey graph shown in \cref{fig:genfareygraph} is induced by the function~$\ell \colon \NN \to \NN$ which maps~$n$ to~$n+1$.
\end{example}

To prove \cref{thm:countablesetsnottypical}, we need two lemmas about grain lines, which are motivated by the proof of \cref{thm:countablesetsnottypical}.
So we prove \cref{thm:countablesetsnottypical} first, giving the statements of the two lemmas where we need them, and then proceed to prove the two lemmas afterwards.

\begin{proof}[Proof of \cref{thm:countablesetsnottypical}]
We have to show that, for every countable set~$\cH$ of infinitely edge-connected graphs, there exists an outerplanar infinitely edge-connected graph~$G$ which contains no subdivision of a graph in~$\cH$.
For this, let~$\cH$ be any countable set of infinitely edge-connected graphs.
By applying \cref{thm:existenceofgrainlines} in each graph in~$\cH$, we find a countably infinite set~$\set{(\cM^{(i)}, \cQ^{(i)}) : i \in \NN}$ of grain lines such that every graph in~$\cH$ contains one of these.
For every~$k \in \NN$, we let~$\ell(k)$ be one greater than the maximum of the lengths of the paths~$Q_j^{(i)}$ with~$0 \leq i,j \leq 2k$.
We remark that~$\ell(0) \geq 1$ and~$\ell(1) \geq 2$, and that the function~$\ell$ is increasing (though not necessarily strictly so).
    
Let~$(\cL, \cP):=(\cL,\cP)(\ell)$ be the grain line induced by~$\ell$, and
let~$G := \hFG\oset{\ell}$ be the generalised halved Farey graph induced by~$\ell$.
In particular, $(\cL, \cP)$ satisfies \cref{GL2'} and \cref{GL3'}, and~$G$ is outerplanar and infinitely edge-connected.
We claim that~$G$ contains no subdivision of a graph in~$\cH$.

Indeed, suppose for a contradiction that some graph~$H$ in~$\cH$ is a topological minor of~$G$.
Then there is a natural number~$i$ such that~$\bigcup \cQ^{(i)}$ is a topological minor of~$G$.
In \cref{subsec:wellstructuredgrainlines}, we will see the crucial fact that 
if a generalised halved Farey graph contains a grain line as a topological minor, then it actually contains the grain line as a subgraph (up to some finite error):

\begin{keylemma} \label{lemma:grainlinetopologicalminorimpliesalmostsubgraph_generalisedfareygraphversion}
    Let~$(\cM,\cQ)$ be a grain line such that a subdivision of~$\bigcup\cQ$ is contained in a generalised halved Farey graph. Then there exists a number~$d$ such that no edge of~$\bigcup\cQ_{\ge d}$ is subdivided.
\end{keylemma}

Hence, by \cref{lemma:grainlinetopologicalminorimpliesalmostsubgraph_generalisedfareygraphversion}, there is a number~$d$ such that~$(\cM, \cQ^{(i)}_{\geq d})$ is a grain line in~$G$.
In \cref{subsec:learnhowtodive}, we will see that while we `dive deeper' into the grain line~$(\cM, \cQ^{(i)}_{\geq d})$, we `dive' with at least the same speed into the grain line~$(\cL,\cP)$:

\begin{keylemma} \label{lemma:simultaneouslydivingdeepintograinlines_generalisedfareygraphversion}
    Let~$(\cM, \cQ)$ be a grain line in the graph~$\bigcup\cP'$ defined by a grain line~$(\cL', \cP')$ satisfying \cref{GL2'} and \cref{GL3'}.
    Then there are a natural number~$q$, a strictly increasing sequence~$(p_k : k \in \NN)$ of natural numbers and a sequence of nested intervals~$([u_k, v_k]_\cM : k \in \NN)$ of~$\cM$ such that, for every~$k \in \NN$, the path~$u_{k+1} P'_{p_k} v_{k+1}$ is a~$\cP'$-segment and a subpath of~$u_k Q_{q+k} v_k$.
\end{keylemma}

Thus, by \cref{lemma:simultaneouslydivingdeepintograinlines_generalisedfareygraphversion} applied to~$(\cL',\cP'):=(\cL,\cP)$, there is a natural number~$q \geq d$, a strictly increasing sequence~$(p_k : k \in \NN)$ of natural numbers and a sequence of nested intervals~$([u_k, v_k]_{\cM^{(i)}} : k \in \NN)$ of~$\cM^{(i)}$ such that, for every~$k \in \NN$, the path~$u_k Q_{q+k}^{(i)} v_k$ contains a~$\cP$-segment~$u_{k+1} P_{p_k} v_{k+1}$.
In particular, by the definition of~$(\cL, \cP)$, the path~$Q_{q + k}^{(i)}$ has length at least~$\ell(p_k) \geq \ell(p_0 + k) \geq \ell(k)$ for every~$k\in\NN$. 
For~$k = \max \bigset{ q, \left\lceil \frac{i}{2} \right\rceil}$, we have~$i \leq 2 \cdot \left\lceil \frac{i}{2} \right\rceil \leq 2k$ and~$q+k \leq 2k$.
Thus, by the definition of~$\ell(k)$, the path~$Q_{q + k}^{(i)}$ is at least one longer than itself, a contradiction.
\end{proof}

\begin{proof}[Proof of \cref{mainthm:anytypicalsetisbig}]
\cref{thm:countablesetsnottypical} implies \cref{mainthm:anytypicalsetisbig}.
\end{proof}

So to complete the proof of \cref{thm:countablesetsnottypical}, it is left to show \cref{lemma:grainlinetopologicalminorimpliesalmostsubgraph_generalisedfareygraphversion} and \cref{lemma:simultaneouslydivingdeepintograinlines_generalisedfareygraphversion}.
We will prove them in \cref{subsec:wellstructuredgrainlines} and in \cref{subsec:learnhowtodive}, respectively.

\subsection{Proof of \texorpdfstring{\cref{lemma:grainlinetopologicalminorimpliesalmostsubgraph_generalisedfareygraphversion}}{Key Lemma 3.4}} \label{subsec:wellstructuredgrainlines}

A grain line~$(\cL, \cP)$ is \emph{well-structured} if, for every~$\cP$-segment~$uP_d v$, we have the inclusion~$V(uP_d v) \cap L \subseteq [u,v]_{\cL}$.
We say that a grain line~$(\cL,\cP)$ is \emph{\free} if~$\bigcup \cP$ is infinitely edge-connected.
Obviously, the following assertions are equivalent:
\[
    \begin{array}{lll}
         \text{(1) }(\cL,\cP)\text{ is \free ;} & \text{(2) }L=V(\bigcup\cP)\text{;} & \text{(3) no vertex of }\bigcup\cP\text{ has degree two.}
    \end{array}
\]
We remark but will not use that free and well-structured grain lines are wild.

The grain line that defines the halved Farey graph in \cref{example:halvedfareyisgrain}, for instance, is both free and well-structured.
More general examples which may be helpful to have in mind while thinking about \free\ and well-structured grain lines are grain lines which satisfy \cref{GL2'} and~\cref{GL3'}.
In fact, the graphs which we constructed in the proof of \cref{thm:countablesetsnottypical} are also defined by such a grain line,\,i.e., they are generalised halved Farey graphs by \cref{lem:genhalvedfareygraphisgrainline}.
We remark that these are not only wild but also wildly presented.

A crucial property of any well-structured grain line\ is that the deletion of any internal vertex destroys the infinite edge-connectivity between its startvertex and endvertex (see \cref{fig:onlyfinitelymanyedgedisjointpaths}):

\begin{figure}[ht]
    \centering
    \includegraphics[width=.5\textwidth]{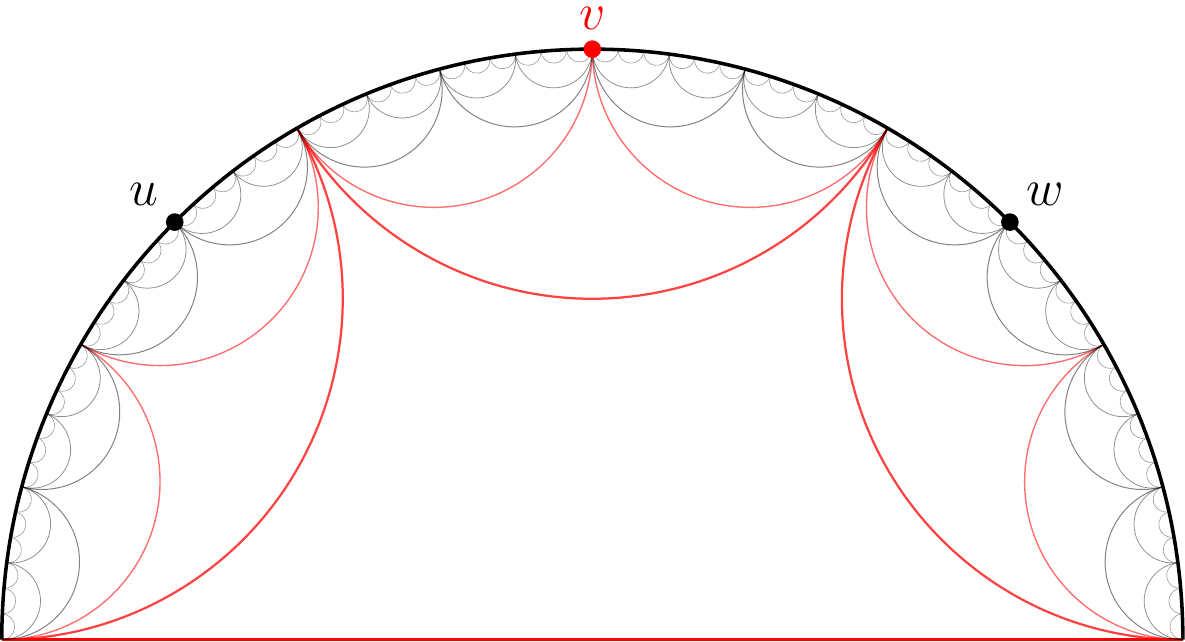}
    \caption{A generalised halved Farey graph in which its vertex~$v$ together with the edges (red) in depth at most the depth of~$v$ separate~$u$ and~$w$.}
    \label{fig:onlyfinitelymanyedgedisjointpaths}
\end{figure}

\begin{lemma} \label{lemma:onlyfinitelymanyedgedisjointpaths}
    Let~$(\cL, \cP)$ be a well-structured grain line and let~$u <_L v <_L w$ be three vertices in~$L$.
    Then the vertex~$v$ together with the set~$F$ of all the finitely many edges in~$\cP$-depth at most the~$\cP$-depth of~$v$ separate~$u$ and~$w$ in~$\bigcup\cP$.
    That is, $v$ together with~$F$ separates all vertices in~$[x,v)_{\cL}$ from all vertices in~$(v,y]_{\cL}$, where~$x$ and~$y$ are the minimum and maximum of~$\cL$, respectively.
\end{lemma}

\begin{proof}
    Let~$d$ be the depth of~$v$ in~$\cP$.
    Since~$(\cL, \cP)$ is well-structured, we have for every path~$P_n$ with~$n > d$
    the two inclusions
    \[
    V(xP_nv - v) \cap L \subseteq [x,v)_\cL\quad\text{ and }\quad V(vP_ny - v) \cap L \subseteq (v,y]_\cL.
    \]
    Hence, every $[x,v)_\cL$--$(v,y]_\cL$~path in~$\bigcup \cP$ which avoids~$v$ has to contain at least one of the finitely many edges in depth at most~$d$ in~$\cP$.
\end{proof}

It follows that every well-structured grain line~$(\cL,\cP)$ imposes the order of~$\cL$ on every grain line in~$\bigcup\cP$:

\begin{lemma} \label{cor:subgrainlinesareintervals}
    Let~$(\cM, \cQ)$ be a grain line in the graph~$\bigcup\cP$ defined by a grain line~$(\cL, \cP)$.
    If~$(\cL, \cP)$ is well-structured, then~$M$ is an interval of~$\cL$ and~$\leq_M$ is either the linear order on~$M$ induced by~$\leq_L$ or the reverse of it.
\end{lemma}

\begin{proof}
    First, we note that~$M$ is a subset of~$L$: By \cref{prop:externallydisjoint}, any vertex of~$\bigcup\cP$ which is not in~$L$ has degree two;
    but every vertex of~$M$ has infinite degree in~$\bigcup \cQ \subseteq \bigcup \cP$ by \cref{prop:Lislimit}. 
    We claim that it suffices to show the following assertion~\cref{Fact:infGhasEitherRayOrInfDefVx}.
    \vspace{.5\baselineskip}
    \begin{equation*}%
    \tag{$\ast$}%
    \hspace{\parindent}\begin{aligned}\label{Fact:infGhasEitherRayOrInfDefVx}%
        \parbox{\textwidth-3\parindent}{\emph{For every three vertices~$u,w\in M$ and~$v\in L$ with~$u\le_L v\le_L w$, the vertex~$v$ is contained in~$M$ and \\ the vertex~$v$ lies in between~$u$ and~$w$ with respect to~$\leq_M$.}}%
    \end{aligned}%
    \end{equation*}%

    \vspace{.5\baselineskip}
    \noindent 
    Let~$x$ and~$y$ be the endvertices of~$(\cM, \cQ)$ named so that~$x <_L y$.
    First, let us deduce~$M = [x,y]_\cL$ from \cref{Fact:infGhasEitherRayOrInfDefVx}.
    By possibly reversing~$\leq_M$, we may assume without loss of generality that also~$x <_M y$.
    For the forward inclusion~$M \subseteq [x,y]_\cL$, let~$v \neq x,y$ be any element of~$M$.
    We show~$x<_L v$ and~$v<_L y$ separately.
    To see that~$x<_L v$, we suppose for a contradiction that~$v<_L x$.
    Then we apply \cref{lemma:onlyfinitelymanyedgedisjointpaths} to the grain line~$(\cL,\cP)$ and~$v <_L x <_L y$ to find a finite set~$F_x$ of edges of~$\bigcup\cP$ such that~$x$ together with~$F_x$ separates~$v$ from~$y$ in the graph~$\bigcup\cP$.
    But~$v\cQ y$ is an infinite system of pairwise edge-disjoint $v$--$y$ paths that avoid~$x$ (since~$x<_M v$), contradicting that~$F_x$ is finite.
    Hence we have~$x <_L v$.
    An analogue argument shows~$v <_L y$.
    Hence, $v \in [x,y]_\cL$.
    Applying \cref{Fact:infGhasEitherRayOrInfDefVx} to~$x$ and~$y$ and each~$v\in [x,y]_\cL$ yields the backward inclusion~$M\supset [x,y]_\cL$.
    Second, we verify that~$\leq_M$ and~$\leq_L$ agree on~$M$.
    For every two vertices~$u, v$ in~$M$ with~$u \le_L v$, it follows from~$M \subseteq [x,y]_\cL$ that~$x \le_L u \le_L v$ and, since~$x \le_M v$, we obtain from~\cref{Fact:infGhasEitherRayOrInfDefVx} that~$x \le_M u \le_M v$.
    
    Now suppose for a contradiction that \cref{Fact:infGhasEitherRayOrInfDefVx} does not hold,\,i.e., suppose that there are three vertices~$u,w \in M$ and~$v \in L$ with~$u \le_L v \le_L w$ such that~$v$ does not lie between~$u$ and~$w$ with respect to~$\leq_M$. Thus, $u <_L v <_L w$.
    Then the sequence~$\oset{uQ_nw : n > d}$ with~$d$ larger than the~$\cQ$-depths of~$u,v$ and~$w$ is an infinite system of pairwise edge-disjoint~$u$--$w$ paths in~$\bigcup \cP$ which avoid~$v$, contradicting \cref{lemma:onlyfinitelymanyedgedisjointpaths} applied to the grain line~$(\cL,\cP)$ and~$u <_L v <_L w$.
\end{proof}

Next, we show that if~$(\cL,\cP)$ is not only well-structured but also free, then~$(\cL,\cP)$ also imposes its freedom on every grain line in~$\bigcup\cP$ (up to some finite and hence negligible error).

\begin{lemma} \label{lemma:grainlineinfreeisalmostfree}
    Let~$(\cM, \cQ)$ be a grain line in the graph~$\bigcup\cP$ defined by a grain line~$(\cL, \cP)$.
    If~$(\cL, \cP)$ is \free\ and well-structured, then only finitely many vertices of~$\bigcup \cQ$ lie outside of~$M$.
    In particular, there is a number~$d$ such that~$(\cM, \cQ_{\geq d})$ is \free .
\end{lemma}

\begin{proof}
    Let~$U$ be the set of all vertices of~$\bigcup\cQ$ that lie outside of~$M$.
    Since~$(\cL,\cP)$ is free, we have~$L=V(\bigcup \cP)$, so~$U\cup M\subset L$.
    By \cref{cor:subgrainlinesareintervals}, there exist~$x,y\in L$ with~$x<_L y$ such that~$[x,y]_{\cL}=M$, and we may assume without loss of generality that~$\le_M$ is induced by~$\le_L$ (possibly after reversing~$\le_M$).
    Hence~$U$ is included in~$[a,x)_{\cL}\cup (y,b]_{\cL}$, where~$a$ and~$b$ denote the minimum and maximum of~$\cL$, respectively.
    
    Let us suppose for a contradiction that~$U$ is infinite.
    Without loss of generality, the intersection of~$U$ with~$[a,x)_{\cL}$ is infinite.
    Since there are only finitely many vertices of~$\bigcup\cQ$ at~$\cQ$-depth~$d$ for each~$d\in\NN$, there exists an infinite subset~$U'\subset U$ such that every two distinct vertices in~$U'$ have distinct~$\cQ$-depths.
    Hence the paths~$uQ_{d(u)}y$ with~$u\in U'$, where~$d(u)$ denotes the~$\cQ$-depth of~$u$, are pairwise edge-disjoint.
    Since~$x$ is the minimum of~$\cM$ and since all~$u$ are distinct from~$x$ by definition of~$U$, all these paths avoid~$x$.
    But by \cref{lemma:onlyfinitelymanyedgedisjointpaths} applied to~$u<_L x<_L y$, the vertex~$x$ together with some finitely many edges separates~$u$ and~$y$ in~$\bigcup\cP$, a contradiction.
\end{proof}

Since subdividing vertices have degree two, we conclude from \cref{lemma:grainlineinfreeisalmostfree} that if a \free\ and well-structured grain line~$(\cL,\cP)$ contains another grain line as a topological minor, then~$(\cL,\cP)$ actually contains it as a subgraph (again up to some finite and hence negligible error):

\begin{lemma} \label{lemma:grainlinetopologicalminorimpliesalmostsubgraph}
    Let~$(\cL, \cP)$ and~$(\cM, \cQ)$ be grain lines such that~$\bigcup \cP$ contains a subdivision of~$\bigcup \cQ$.
    If~$(\cL, \cP)$ is \free\  and well-structured, then there is a number~$d$ such that~$(\cM, \cQ_{\geq d})$ is a grain line in~$\bigcup \cP$. \qed
\end{lemma}

\begin{proof}[Proof of \cref{lemma:grainlinetopologicalminorimpliesalmostsubgraph_generalisedfareygraphversion}]
    Since a generalised halved Farey graph is defined by a free and well-structured grain line, this follows directly from \cref{lemma:grainlinetopologicalminorimpliesalmostsubgraph}. 
\end{proof}

\subsection{Proof of \texorpdfstring{\cref{lemma:simultaneouslydivingdeepintograinlines_generalisedfareygraphversion}}{Key Lemma 3.5}} \label{subsec:learnhowtodive}

We need the following three lemmas.

\begin{lemma} \label{lemma:firstdive}
    Let~$(\cL, \cP)$ be a grain line and~$P$ a path in~$\bigcup \cP$.
    If~$P$ contains at least one edge in~$\cP$-depth greater than the~$\cP$-depth of its endvertices, then it already contains a~$\cP$-segment in depth equal to the maximum of the~$\cP$-depths of the edges of~$P$.
\end{lemma}

\begin{proof}
Let~$d$ be the maximum of the depths of the edges of~$P$ in~$\cP$.
Let~$uPv$ be a subpath of~$P$ that is a~$(\bigcup\cP_{<d})$-path;
this subpath exists because an edge of~$P$ has~$\cP$-depth~$d$ and the endvertices of~$P$ are contained in~$\bigcup\cP_{<d}$ by assumption.
It follows directly from the choice of~$uPv$ that all of its edges have~$\cP$-depth at least~$d$, and hence~$uPv$ actually is a subpath of~$P_d$. 
So both~$u$ and~$v$ lie in~$L$ by~\cref{prop:Lislimit}, and thus in~$L_{<d}$.
Hence, $uPv$ is the desired~$\cP$-segment.
\end{proof}

\begin{lemma}\label{lem:nestedintervals}
    Let~$(\cL,\cP)$ be a grain line and let~$P$ be a~$u$--$v$ path in~$\bigcup \cP$ containing an edge~$e = ww'$ such that~$u \leq_L w <_L w' \leq_L v$ and such that the~$\cP$-depth~$d$ of~$e$ is greater than the maximum of the~$\cP$-depths of~$u$ and~$v$.
    Let~$Q$ be the subpath of~$P$ which contains~$e$ and is a~$(\bigcup \cP_{<d})$-path.
    If~$(\cL, \cP)$ is well-structured, then $V(Q) \cap L \subseteq [u,v]_\cL$.
\end{lemma}

\begin{proof}
    We denote the minimum and maximum of~$\cL$ by~$a$ and~$b$, respectively.
    It suffices to show that~$Q$ avoids both~$[a,u)_{\cL}$ and~$(v,b]_{\cL}$.
    By symmetry, it suffices to show that~$Q$ avoids~$[a,u)_{\cL}$.
    Let~$F$ denote the set of edges of~$\bigcup\cP_{<d}$.
    By \cref{lemma:onlyfinitelymanyedgedisjointpaths}, the vertex~$u$ separates~$[a,u)_{\cL}$ from~$(u,b]_{\cL}$ in the graph~$(\bigcup\cP)-F$.
    The path~$Q$ contains the vertex~$w'\in (u,b]_{\cL}$; it avoids~$F$ because~$Q$ is a~$(\bigcup\cP_{<d})$-path; and it does not contain~$u$ as an internal vertex because~$Q$ is a subpath of the~$u$--$v$ path~$P$.
    Therefore, $Q$ avoids~$[a,u)_{\cL}$.
\end{proof}

\begin{lemma} \label{lemma:simultaneouslydivingdeepintograinlines}
    Let~$(\cM, \cQ)$ be a grain line in the graph~$\bigcup\cP$ defined by a free, well-structured and wildly presented grain line~$(\cL, \cP)$.
    Then there are a natural number~$q$, a strictly increasing sequence~$(p_k : k \in \NN)$ of natural numbers and a sequence of nested intervals~$([u_k, v_k]_\cM : k \in \NN)$ of~$\cM$ such that, for every~$k \in \NN$, the path~$u_{k+1} P_{p_k} v_{k+1}$ is a~$\cP$-segment and it is a subpath of~$u_k Q_{q+k} v_k$.
\end{lemma}

\begin{proof}
    By \cref{lemma:grainlineinfreeisalmostfree}, we may assume by choosing~$q$ large enough that, without loss of generality, $(\cM, \cQ)$ is free.
    Let~$u_0, v_0 \in \bigcup \cP$ be the two vertices for which~$(\cM, \cQ)$ is a~$u_0$--$v_0$ grain line.
    By \cref{cor:subgrainlinesareintervals} and possibly interchanging~$u_0$ and~$v_0$, we may assume without loss of generality that~$\leq_L$ and~$\leq_M$ agree on~$M\subset L=V(\bigcup\cP)$.
    Moreover, the lemma also ensures that~$M$ is an interval of~$\cL$.
    Let~$p$ be the maximum of the~$\cP$-depths of~$u_0$ and~$v_0$.
    Since there are only finitely many edges in~$\cP$-depth at most~$p$, we can pick~$q$ large enough so that~$Q_q$ contains an edge of~$\bigcup\cP$ in~$\cP$-depth greater than~$p$.
    Then by \cref{lemma:firstdive}, there is a number~$p_0 > p$ such that~$Q_q$ contains a~$\cP$-segment~$u_1P_{p_0}v_1$.
    Since~$(\cM, \cQ)$ is free, $u_1$ and~$v_1$ are in~$M$.
    Thus, $[u_1,v_1]_{\cM} \subseteq [u_0,v_0]_{\cM}$ because~$u_0$ and~$v_0$ are the minimum and maximum of~$\cM$.

    Since~$(\cL,\cP)$ is wildly presented, there is an internal vertex~$w_1$ of~$u_1P_{p_0}v_1$ which is contained in~$(u_1,v_1)_\cL$.
    Since~$(\cM, \cQ)$ is free and~$\leq_L$ and~$\leq_M$ agree on~$M$, the vertex~$w_1$ is contained in~$M$ and~$w_1 \in (u_1,v_1)_\cM$.
    Thus by \cref{prop:Lislimit}, $Q_{q+1}$ contains all three vertices~$u_1,v_1,w_1$ and, by \cref{prop:samelinearorder}, $w_1$ is an internal vertex of~$u_1Q_{q+1}v_1$.
    Since~$u_1 P_{p_0}v_1$ is a~$\cP$-segment, the~$\cP$-depth of the internal vertex~$w_1$ of~$u_1 P_{p_0}v_1$ is~$p_0$.
    Since~$Q_q$ already contains the two edges of~$P_{p_0}$ that are incident with~$w_1$, and since~$w_1$ lies on no path~$P_p$ with~$p<p_0$, the path~$u_1Q_{q+1}v_1$ contains an edge~$e_{q+1}$ incident with~$w_1$ in~$\cP$-depth at least~$p_0+1$.
    Let~$Q'_{q+1}$ be the subpath of~$u_1Q_{q+1}v_1$ which contains~$e_{q+1}$ and is a~$(\bigcup \cP_{\leq p_0})$-path.
    By \cref{lemma:firstdive}, $Q'_{q+1} \subseteq u_1Q_{q+1}v_1$ contains a~$\cP$-segment~$u_2P_{p_1} v_2$ with~$p_1 > p_0$.
    By \cref{lem:nestedintervals} and since~$(\cL, \cP)$ is free, $V(Q'_{q+1}) \subseteq [u_1,v_1]_\cM$.
    Since~$u_2P_{p_1} v_2$ is a subpath of~$Q'_{q+1}$, we have~$[u_2,v_2]_{\cM}\subseteq [u_1,v_1]_{\cM}$.
    
    Next, we repeat the previous step where we replace~$u_1P_{p_0}v_1$, $q$ and~$p_0$ with~$u_2P_{p_1} v_2$, $q+1$ and~$p_1$, and iterate in this way to conclude the proof.
\end{proof}

\begin{proof}[Proof of \cref{lemma:simultaneouslydivingdeepintograinlines_generalisedfareygraphversion}]
    Since a grain line which satisfies \cref{GL2'} and \cref{GL3'} is free, well-structured and wildly presented, this follows directly from \cref{lemma:simultaneouslydivingdeepintograinlines}.
\end{proof}

\subsection{Excluding \texorpdfstring{$k$}{k}-bounded Farey graph minors}
In~\cite{kurkofka2020fareydeterminedbyconnectivity}, it was shown that the Farey graph is uniquely determined by its connectivity, as follows.
A~\emph{\Pigraph} is an infinitely edge-connected graph such that no two of its vertices are linked by infinitely many pairwise internally disjoint paths.
A~\emph{$\kappa$-bounded} minor, for a cardinal~$\kappa$, is a minor with branch sets of size less than~$\kappa$.
A~\Pigraph\ is \emph{$\kappa$-typical} if it occurs as a~$\kappa$-bounded minor in every \Pigraph . 
Note that any two~$\kappa$-typical \Pigraph s are~$\kappa$-bounded minors of each other; we call such graphs~$\kappa$-boundedly minor-equivalent.

\begin{theorem}\cite{kurkofka2020fareydeterminedbyconnectivity}\label{ConnectivityChar}
    Up to~$\aleph_0$-bounded minor-equivalence, the Farey graph is the unique~$\aleph_0$-typical \Pigraph .
\end{theorem}

A~referee of~\cite{kurkofka2020fareydeterminedbyconnectivity} asked whether~$\aleph_0$ is best possible for the above theorem; that is, can~$\aleph_0$ be replaced in~\cref{ConnectivityChar} with some~$\kappa<\aleph_0$ so that the theorem's statement remains true?
Using generalised halved Farey graphs, we can show that the answer is in the negative:

\begin{lemma}
For every~$k\in\NN$ there exists a \Pigraph\ that does not contain the Farey graph as a~$k$-bounded minor and that is a generalised halved Farey graph.
\end{lemma}
\begin{proof}
Let~$k\in\NN$ be given.
On the one hand, as the Farey graph contains a triangle, every graph that contains the Farey graph as a~$k$-bounded minor must contain a cycle of length at most~$3k$.
On the other hand, every generalised halved Farey graph clearly is a \Pigraph .
So it suffices to find a generalised halved Farey graph of girth at least~$3k+1$.
For example, we could take~$\hFG(\ell)$ for~$\ell\colon\NN\to\NN$ given by~$\ell(0):=1$ and~$\ell(n):=3k$ when~$n\ge 1$.
\end{proof}

\section{Typical graph with regard to the immersion relation}\label{sec:StrongImmersionProof}

\noindent In this section, we prove \cref{mainthm:halvedFareygraphistypical}: we show that the halved Farey graph is immersed in every infinitely edge-connected graph.
Before we do this, however, we take a step back to verify that the other two obvious candidates, namely~$\Kaleph$ and the Farey graph, are not immersed in every infinitely edge-connected graph.
For readers who are familiar with~\cite{kurkofka2020infinitelyedgecongraphcontaingfareygraphorT}, we remark that the second infinitely edge-connected graph mentioned in \cref{thm:originalFarey} is immersion-equivalent to~$\Kaleph$; in particular, it is not immersed in every infinitely edge-connected graph either.

\begin{theorem} \label{thm:usualcandidatesandimmersion}
    \text{}
    \begin{enumerate}
        \item \label{item:infinitecomplete} $\Kaleph$ is not immersed in the Farey graph.
        \item \label{item:Farey} Neither~$\Kaleph$ nor the Farey graph is immersed in the halved Farey graph.
    \end{enumerate}
\end{theorem}

\begin{proof}
    \cref{item:infinitecomplete}: Suppose for a contradiction that there is an immersion of~$\Kaleph$ in the Farey graph~$\FG$.
    Let~$U$ be the set of branch vertices of this immersion.
    By the definition of the Farey graph, we may write~$\FG$ as the union of two halved Farey graphs which intersect in the Farey graph of order~$0$.
    Hence, the two linear orders on the vertex sets of these two halved Farey graphs introduced in \cref{example:halvedfareyisgrain} induce a cyclic order on the vertex set of~$\FG$.
    Now, there are two branch vertices~$u,v \in U$ such that the open intervals~$(u,v)$ and~$(v,u)$ with respect to this cyclic order intersect~$U$ non-emptily.
    Since~$\Kaleph$ is infinitely edge-connected even after deleting the two vertices~$u$ and~$v$, the immersion yields infinitely many edges between~$(u,v)$ and~$(v,u)$ in~$\FG$.
    This contradicts the fact that, by the construction of~$\FG$, there are only finitely many of these.

    \cref{item:Farey}: Since the halved Farey graph~$\hFG$ is a subgraph of the Farey graph, $\Kaleph$ cannot be immersed in~$\hFG$ by~\cref{item:infinitecomplete}.
    Suppose for a contradiction that there is an immersion of the Farey graph~$\FG$ in the halved Farey graph~$\hFG$.
    Let~$U$ be the set of branch vertices of this immersion.
    Let~$x,y$ be the two vertices in the halved Farey graph of order~$0$ in the halved Farey graph~$\hFG$.
    There is a branch vertex~$w \in U$ such that the two intervals~$[x,w)$ and~$(w,y]$ with respect to the linear order on~$V(\hFG)$ introduced in \cref{example:halvedfareyisgrain} intersect~$U$ non-emptily.
    Since the Farey graph~$\FG$ is infinitely edge-connected even after deleting the single vertex~$w$, the immersion yields infinitely many edges between~$[x,w)$ and~$(w,y]$ in~$\hFG$.
    This contradicts the fact that, by the construction of~$\hFG$, there are only finitely many of those.
\end{proof}

Now we tend to the proof of \cref{mainthm:halvedFareygraphistypical}.

\subsection{Overview of the proof of Theorem 1}

Our aim for the remainder of this paper is to show that every infinitely edge-connected graph contains the halved Farey graph as an immersion minor.
For this, it obviously suffices to consider only infinitely edge-connected graphs with no~$\Kaleph$-immersion-minors.

As our first step, we will transfer the notion of cutvertices from vertex-connectivity to infinite edge-connectivity by introducing `\ccvs ', and we will introduce what could be considered an analogue of the block-cutvertex theorem for \ccvs\ and blocks of infinite edge-connectivity (\cref{lemma:nestedsetseparatingeverything}).
Unlike the tree-structure of a graph imposed by its cutvertices, the structure imposed by its \ccvs\ will be tree-like but not in general a tree-decomposition.
Indeed, the tree-like structure that we will obtain from the \ccvs\ of a graph can exhibit~$(\omega+1)$-chains, and for some graphs like the halved Farey graph, their tree-like structure can even be order-isomorphic to~$\QQ$.

As our second step, we will show that an infinitely edge-connected graph contains the halved Farey graph as an immersion minor if the tree-like structure given by its \ccvs\ is `wild' in that it exhibits an interval which is order-isomorphic to~$\QQ$ (\cref{cor:Qorderedgivesfarey}).

So as our third and final step, we will deal with the case in which the tree-like structure given by the \ccvs\ is not wild. 
Roughly, we will employ the \ccvs\ to find infinitely many infinitely edge-connected subgraphs which are almost vertex-disjoint and which have no \ccvs\ themselves.
Then we will choose an arbitrary vertex in each subgraph that lies in no other of these subgraphs, and we will link the chosen vertices up with paths to obtain a~$\Kaleph$-immersion-minor.
As this contradicts our initial assumption that no~$\Kaleph$-immersion-minor is present, the proof will be concluded.

\subsection{Compound separations}\label{sec:introofcsep}

Recall that a \emph{separation of a set}~$V$ is a set~$\{A,B\}$ such that~$A\cup B=V$. 
We call~$A$ and~$B$ the \emph{sides} of this separation.
The separation is \emph{proper} if~$A\setminus B$ and~$B\setminus A$ are non-empty.
Now let~$G$ be any infinite graph.
A~\emph{\csep } of~$G$ is a proper separation~$\{A,B\}$ of~$V(G)$ such that the \emph{\cseparator}~$A\cap B$ is finite and~$G$ has only finitely many edges between~$A\setminus B$ and~$B\setminus A$.
Then the cardinality of the separator~$A\cap B$ is the \emph{order} of the \csep~$\{A,B\}$.
If the separator of a \csep~$\{A,B\}$ is a singleton~$\{u\}$, then we also refer to the vertex~$u$ as the \emph{separator} of~$\{A,B\}$, and we say that~$\{A,B\}$ is \emph{\unit}.
A~vertex~$u$ is a \emph{\ccv} of~$G$ if there exists a \csep\ of~$G$ with separator~$u$.
A \csep~$\{A,B\}$ of~$G$ \emph{separates} two vertices~$u,v$ of~$G$ if~$u$ is contained in~$A\setminus B$ and~$v$ is contained in~$B\setminus A$, or vice versa.
It separates~$u$ and~$v$ \emph{minimally} if no \csep~$\{C,D\}$ with~$C\cap D\subsetneq A\cap B$ separates~$u$ and~$v$ in~$G$.
For example, if~$\{A,B\}$ is a \unit\ \csep\ of an infinitely edge-connected graph with~$u\in A\setminus B$ and~$v \in B\setminus A$, then~$\{A,B\}$ minimally separates~$u$ and~$v$ in that graph.
We say that two vertices~$u$ and~$v$ of~$G$ are \emph{$k$-\ccon } in~$G$ for a natural number~$k$ if no \csep\ of~$G$ of order less than~$k$ separates~$u$ and~$v$.
If every two vertices of~$G$ are~$k$-\ccon , then~$G$ itself is \emph{$k$-\ccon }.

If an infinite graph~$G$ is~$k$-\ccon\ for every~$k \in \NN$, then~$G$ is infinitely vertex-connected; in particular, we greedily find~$\Kaleph$ as a topological minor in~$G$, so the halved Farey graph is immersed in~$G$:

\begin{observation}\label{fact:findcsep}
    Let~$G$ be an infinite graph.
    If the halved Farey graph is not immersed in~$G$, then there is a pair of vertices of~$G$ which is separated by some \csep\ of~$G$.
    \qed
\end{observation}

Every minimally separating \csep\ of an infinitely edge-connected graph~$G$ can be used to split~$G$ into two infinitely edge-connected subgraphs:

\begin{keylemma} \label{lemma:splittinglemma}
    Let~$G$ be an infinitely edge-connected graph and~$\set{A,B}$ a \csep\ of~$G$ that minimally separates two vertices of~$G$.
    Then~$G[A]$ and~$G[B]$ are infinitely edge-connected.
\end{keylemma}

We prepare the proof of \cref{lemma:splittinglemma} with the following lemma:

\begin{lemma} \label{lemma:connectivitylemmafortightmixseps}
    Let~$G$ be an infinitely edge-connected graph and~$\set{A,B}$ a \csep\ of~$G$ that minimally separates two vertices~$u,v \in G$.
    If~$u \in A$, then there is a system of pairwise edge-disjoint~$u$--$(A \cap B)$ paths in~$G[A]$ such that, for every vertex~$w \in A \cap B$, infinitely many of these paths end in~$w$.
\end{lemma}

\begin{proof}
    It suffices to find for every vertex~$w \in A \cap B$ a system~$\cP$ of infinitely many pairwise edge-disjoint~$u$--$w$ paths in~$G[A]$ which are internally disjoint from~$A \cap B$.
    For this, let~$w \in A \cap B$ be any vertex.
    Set~$S := (A \cap B) \setminus \set{w}$ and~$H := G - S$.
    Any finite cut~$E(C,D)$ of~$H$ which separates~$u$ and~$v$ induces a \csep~$\set{C \cup S, D \cup S}$ of~$G$ with separator~$S \subsetneq A \cap B$ which separates~$u$ and~$v$.
    Since the latter does not exist by assumption, $u$ and~$v$ cannot be separated by finitely many edges in~$H$.
    Hence, there are infinitely many pairwise edge-disjoint~$u$--$v$ paths in~$G$ that avoid~$S$.
    Since~$\set{A,B}$ is a \csep\ that separates~$u$ and~$v$, all but finitely many of these paths have to meet~$w$.
    Thus, the~$u$--$w$ subpaths of these paths form the desired path system.
\end{proof}

\begin{proof}[Proof of \cref{lemma:splittinglemma}]
    By symmetry, it is enough to prove that~$G[A]$ is infinitely edge-connected.
    Let~$u$ and~$v$ be two vertices of~$G$ which are minimally separated by a \csep~$\{A,B\}$, named so that~$u\in A\setminus B$.
    To show that~$G[A]$ is infinitely edge-connected, it suffices to find infinitely many pairwise edge-disjoint~$u$--$a$ paths in~$G[A]$ for every vertex~$a\in A$ other than~$u$. For this, let any such vertex~$a$ be given.
    Since~$G$ is infinitely edge-connected, there is an infinite system~$\cQ'$ of pairwise edge-disjoint~$a$--$v$ paths in~$G$.
    As~$\set{A, B}$ is a \csep\ of~$G$, it follows from the pigeonhole principle that there is a vertex~$w \in A \cap B$ for which there are infinitely many paths in~$\cQ'$ whose first vertex in~$B$ is~$w$.
    We denote the system of the~$a$--$w$~subpaths of these paths by~$\cQ$.
    Note that all paths in~$\cQ$ are included in~$G[A]$.
    By \cref{lemma:connectivitylemmafortightmixseps}, there is an infinite system~$\cP$ of pairwise edge-disjoint~$u$--$w$ paths in~$G[A]$.
    Since~$\cQ$ and~$\cP$ both consist of paths in~$G[A]$, we can greedily combine these two path systems to obtain infinitely many pairwise edge-disjoint~$u$--$a$ paths in~$G[A]$.
    Thus, $G[A]$ is infinitely edge-connected.
\end{proof}

\subsection{Faithful nested sets of \unit\ \csep s} \label{sec:constructnestedset}

In this section, we prove the following key lemma which could be viewed as an analogue of the block-cutvertex theorem for infinite edge-connectivity. Missing definitions follow.

\begin{keylemma}\label{lemma:nestedsetseparatingeverything} 
    Let~$G$ be an infinitely edge-connected graph.
    If~$\Kaleph$ is not immersed in~$G$, then there is a nested set of \unit\ \csep s of~$G$ which is faithful to~$G$.
\end{keylemma}

First, we make this statement precise.
Let~$G$ be an infinitely edge-connected graph.
For a separation~$\set{A,B}$ of a set~$V$, we recall that the ordered pairs~$(A,B)$ and~$(B,A)$ are its \emph{orientations}, and~$(A,B)$ and~$(B,A)$ are \emph{oriented separations} of~$V$.
For two oriented separations~$(A,B)$ and~$(C,D)$, we have that~$(A,B) \leq (C,D)$ if~$A \subseteq C$ and~$B \supseteq D$.
Two separations are \emph{nested} if they have comparable orientations.
Moreover, a set of separations is \emph{nested} if every two separations in it are nested.
By choosing precisely one orientation of each unoriented separation in a set~$S$ of unoriented separations, we obtain an \emph{orientation}~$\sigma$ of~$S$. 
The orientation~$\sigma$ is a \emph{star} of separations if for every two distinct~$(A,B), (C,D) \in \sigma$ we have~$(A,B) \leq (D,C)$.

Let~$w$ be a \ccv\ of~$G$.
A nested set~$N$ of \unit\ \csep s of~$G$ is \emph{faithful} to~$w$ in~$G$ if
\begin{itemize}
    \item the separators of all separations in~$N$ are equal to~$w$,
    \item there is some orientation of~$N$ which is a star of separations,
    \item and for every pair of vertices which are separated by some \csep\ with separator~$w$ there is some \csep\ in~$N$ separating them.
\end{itemize}
A nested set~$N$ of \unit\ \csep s of~$G$ is \emph{faithful} to~$G$ if for every \ccv~$w$ of~$G$ the subset~$N_w$ of~$N$ containing all \csep s of~$N$ with separator~$w$ is faithful to~$w$ in~$G$.

If we can find for each \ccv~$w$ of~$G$ a nested set of \csep s of~$G$ which is faithful to~$w$ in~$G$, then we can combine these to obtain a nested set which is faithful to~$G$:

\begin{lemma} \label{lemma:distinctseparatorsarenested}
    Let~$G$ be an infinitely edge-connected graph.
    If~$\set{A, B}$ and~$\set{C, D}$ are two \unit\ \csep s of~$G$ with distinct separators~$u$ and~$v$, respectively, such that~$u \in C$ and~$v \in B$, then~$(A,B) \leq (C,D)$ and, in particular, $\set{A, B}$ and~$\set{C, D}$ are nested.
\end{lemma}

\begin{proof}
    The induced subgraph~$G[D]$ shares with~$B$ at least the vertex~$v$ by assumption, but it does not contain the vertex~$u$ because~$u$ is contained in~$C\setminus\{v\}=V(G)\setminus D$. Since~$G[D]$ is infinitely edge-connected by \cref{lemma:splittinglemma} and~$\{A,B\}$ is a \csep\ with separator~$u$, it follows that~$G[D]$ is included in~$G[B]-u$.
    By symmetry, $G[A]$ is included in~$G[C]-v$.
    Hence~$(A,B)\leq (C,D)$.
\end{proof}

Hence to find a nested set of \csep s which is faithful to~$G$, it suffices to find for every \ccv~$w$ of~$G$ a nested set of \csep s which is faithful to~$w$ in~$G$.

The \emph{edge-blocks} of a graph~$G$ are the classes of the equivalence relation `not separable by finitely many edges' on the vertex set of~$G$.
So a vertex~$w$ of~$G$ is a \ccv\ if and only if~$G-w$ is not infinitely edge-connected, if and only if~$G-w$ has at least two edge-blocks.

\begin{lemma}\label{fact:nestedsetforasetS}
    Let~$G$ be an infinitely edge-connected graph and~$w$ be a \ccv\ of~$G$.
    If~$G-w$ has only finitely many edge-blocks, then there is a nested set~$N_w$ which is faithful to~$w$ in~$G$.
\end{lemma}

\begin{proof}
    Every edge-block~$X$ of~$G-w$ induces the separation~$\set{X \cup \set{w}, V(G) \setminus X}$ of~$V(G)$ with separator~$w$, and the set~$N_w$ containing all of these separations is nested.
    Moreover, if we orient the separation induced by an edge-block~$X$ of~$G-w$ as~$\oset{X \cup \set{w}, V(G) \setminus X}$, then the orientation of~$N_w$ containing all of these oriented separations is a star of separations.
    Since every two vertices of~$G$ that are separated by some \csep\ with separator~$w$ are contained in distinct edge-blocks of~$G-w$, it is enough to show that every separation in~$N_w$ is a \csep .
    For this, let~$X$ be an edge-block of~$G-w$. 
    Then~$X$ sends only finitely many edges to each other edge-block of~$G-w$. Since~$G-w$ has only finitely many edge-blocks, only finitely many edges run between~$X$ and~$G-w-X$. Hence~$\set{X \cup \set{w}, V(G) \setminus X}$ is a \csep .
\end{proof}
 
For the proof of \cref{lemma:nestedsetseparatingeverything}, we only need one more ingredient:

\begin{lemma} \label{lemma:dualityforfiniteTGS}
    Let~$G$ be an infinitely edge-connected graph and~$X$ a finite set of vertices of~$G$.
    Then the following assertions are complementary:
    \begin{enumerate}
        \item \label{item:Kaleph} There is an immersion of~$\Kaleph$ in~$G$ with at most one branch vertex in every edge-block of~$G-X$;
        \item \label{item:fin} $G-X$ has only finitely many edge-blocks.
    \end{enumerate}
\end{lemma}

To prove \cref{lemma:dualityforfiniteTGS}, we need the following tool, and the notion of tree-cut decompositions by Wollan~\cite{Wollan}.
Recall that a \emph{near-partition} of a set~$V$ is a family of pairwise disjoint subsets~$X_\xi\subset V$, possibly empty, such that~$\bigcup_\xi X_\xi=V$.
Let~$G$ be a graph, $T$ a tree, and let~$\cX=(X_t)_{t\in T}$ be a family of vertex sets~$X_t\subset V(G)$ indexed by the nodes~$t$ of~$T$.
The pair~$(T,\cX)$ is called a \emph{\tcd } of~$G$ if~$\cX$ is a near-partition of~$V(G)$.
The vertex sets~$X_t$ are the \emph{parts} of the \tcd~$(T,\cX)$.
We say that~$(T, \cX)$ is a \tcd\ \emph{into} these parts.
Whenever a \tcd\ is introduced as~$(T,\cX)$, we tacitly assume that~$\cX=(X_t)_{t\in T}$.
If~$(T,\cX)$ is a \tcd , then we associate with every edge~$t_1 t_2 \in T$ its \emph{adhesion set}~$E_G(\,\bigcup_{t\in T_1}X_t\,,\,\bigcup_{t\in T_2}X_t\,)$ where~$T_1$ and~$T_2$ are the two components of~$T-t_1 t_2$ with~$t_1\in T_1$ and~$t_2\in T_2$.
A \tcd\ has \emph{finite adhesion} if all its adhesion sets are finite.

\begin{theorem}[\,\cite{kurkofka2020infinitelyedgecongraphcontaingfareygraphorT}*{Theorem~5.1}\,]\label{thm:tcdintoinftyblocks}
    Every connected graph has a \tcd\ of finite adhesion into its edge-blocks.
\end{theorem}

We remark that, since edge-blocks are non-empty, the parts of the \tcd\ in \cref{thm:tcdintoinftyblocks} form a partition of the vertex set.
Hence, every adhesion set of this \tcd\ is a cut.

\begin{proof}[Proof of \cref{lemma:dualityforfiniteTGS}]
    Assertions~\cref{item:Kaleph} and \cref{item:fin} clearly exclude each other, so it suffices to prove~$\neg$\cref{item:fin}$\to$\cref{item:Kaleph}.
    
    We assume that~$G-X$ has infinitely many edge-blocks.  
    Then we find an infinite set~$\set{U_n : n \in \NN}$ of subsets of~$V(G) \setminus X$ such that all~$\{U_n \cup X,V(G)\setminus U_n \}$ are \csep s of~$G$ (with separator~$X$), as follows.
    By \cref{thm:tcdintoinftyblocks}, $G-X$ admits a \tcd~$(T, \cX)$ of finite adhesion into its edge-blocks.
    By König's Infinity Lemma (see e.g.~\cite{DiestelBook}*{Lemma 8.1.2}), the tree~$T$ contains a vertex of infinite degree or a ray.
    If~$t_0,t_1,\ldots$ are infinitely many neighbours of some vertex~$t \in T$, then for every~$n\in\NN$ we let~$U_n$ be the side of the finite cut of~$G-X$ induced by the edge~$t_n t \in T$ that does not include~$X_t$.
    Otherwise, we find a ray~$t_0 t_1\ldots$ in~$T$.
    We denote the finite cut of~$G-X$ induced by the edge~$t_nt_{n+1} \in T$ by~$\set{A_n, B_n}$ so that~$A_n$ includes~$X_{t_n}$. 
    We set~$U_0 := A_0$ and, for every natural number~$n \geq 1$, we let~$U_n := A_n \cap B_{n-1}$.
    It is straightforward to see that the resulting partition~$\{U_n: n\in\NN\}$ is as desired in either case.
    
    For every~$n \in \NN$, we fix a vertex~$u_n \in U_n$ and an infinite system~$\cP_n$ of pairwise edge-disjoint~$u_n$--$X$ paths in~$G$.
    Since~$\{U_n \cup X,V(G)\setminus U_n \}$ is a \csep\ of~$G$, all but finitely many of these paths are contained in~$G[U_n \cup X]$; thus, we may assume without loss of generality that all of them are contained in~$G[U_n \cup X]$.
    By the pigeonhole principle, there is a vertex~$x$ in the finite vertex set~$X$ such that there is an infinite set~$M$ of natural numbers~$n$ for which infinitely many of the paths in~$\cP_n$ end in~$x$.
    Now we find an immersion of~$\Kaleph$ in~$G$ with $\set{u_n : n \in M}$ as its set of branch vertices as follows:
    we greedily connect every two branch vertices~$u_i$ and~$u_j$ with a path that avoids all other branch vertices and the already chosen finitely many edges using the path systems~$\cP_i$ and~$\cP_j$.
\end{proof}

\begin{proof}[Proof of \cref{lemma:nestedsetseparatingeverything}]
    Let~$G$ be an infinitely edge-connected graph and suppose that~$\Kaleph$ is not immersed in~$G$.
    For every \ccv~$w$ of~$G$, the graph~$G-w$ has only finitely many edge-blocks by \cref{lemma:dualityforfiniteTGS}, and so
    there is a nested set~$N_w$ which is faithful to~$w$ in~$G$ by \cref{fact:nestedsetforasetS}.
    The union~$N$ of all the sets~$N_w$ is nested by \cref{lemma:distinctseparatorsarenested}, and~$N$ is faithful to~$G$ by construction.
    We remark that, if~$G$ has no \ccv , then~$N=\emptyset$ is faithful to~$G$.
\end{proof}

\subsection{Finding the halved Farey graph in wild structures} \label{sec:grainlinesandstrongimm}

Our next aim is to prove this:

\begin{keylemma} \label{cor:Qorderedgivesfarey}
    Let~$G$ be an infinitely edge-connected graph and~$S$ a set of oriented \unit\ \csep s with pairwise distinct separators. If~$(S,\leq)$ is order-isomorphic to~$\QQ$, then the halved Farey graph is immersed in~$G$.\qed
\end{keylemma}

We have split the proof into two halves which are represented by the following two lemmas:

\begin{lemma}\label{lemma:fareygraphinwildgrainline}
    The halved Farey graph is immersed in any graph that is defined by a wildly presented grain~line.
\end{lemma}

\begin{proof}
    Let~$(\cL, \cP)$ be a wildly presented grain line.
    We define the set~$U$ of branch vertices as follows.
    First, let~$U_0 := \{x,y\}$, where~$x$ and~$y$ are the two endvertices of the grain line.
    Assume that we have already chosen~$U_0, \dots, U_n$ for some natural number~$n$ such that for every~$m \in \set{1, \dots, n}$
    \begin{itemize}
        \item  all vertices in~$U_m \setminus U_{m-1}$ are in~$\cP$-depth~$m$, and
        \item  for every~$u \in U_{m-1}$ with successor~$u'$ in~$\partordrest{\cL}{U_{m-1}}$, there is precisely one~$v \in U_m$ such that~$u <_L v <_L u'$.
    \end{itemize}
    We order~$U_n$ according to~$\leq_L$ as~$u_0 <_L u_1 <_L \cdots <_L u_\ell$.
    Now we obtain~$U_{n+1}$ from~$U_n$ by adding exactly one vertex~$v_k \in L$ for every~$k \in \set{1, \dots, \ell}$ such that~$v_k$ is an internal vertex of the subpath~$u_{k-1}P_{n+1}u_k$ and~$u_{k-1} <_L v_k <_L u_k$:
    Let~$u_k'$ be the~$\leq_L$-minimal vertex in~$L$ in~$\cP$-depth at most~$n$ with~$u_{k-1} <_L u_k' \leq_L u_k$.
    Since~$(\cL, \cP)$ is wildly presented, there is an internal vertex~$v_k$ of~$u_{k-1} P_{n+1} u_k'$ with~$u_{k-1} <_L v_k <_L u_k'$.
    In particular, $v_k <_L u_k$~and, by the choice of~$u_k'$, the vertex~$v_k$ is in~$\cP$-depth~$n+1$.
    
    Then we set~$U := \bigcup_{n \in \NN} U_n$.
    By construction we have that, for every~$n \in \NN$, the path~$P_{n}$ induces the same linear order on~$U_n$ as~$\leq_L$ and, since~$U \setminus U_n$ consists only of vertices in~$\cP$-depth greater than~$n$, the path~$P_n$ is also disjoint from~$U \setminus U_n$.
    Thus, we have found an immersion of the halved Farey graph in~$\bigcup \cP$ with~$U$ as its set of branch vertices.
\end{proof}

\begin{lemma} \label{lemma:Qorderedsepgivewildgrainline}
    Let~$G$ be an infinitely edge-connected graph.
    Let~$S$ be a set of oriented \unit\ \csep s of~$G$ with pairwise distinct separators. If~$(S, \leq)$ is order-isomorphic to~$\QQ \cap [0,1]$, then there exists a wildly presented 
    grain line~$(\cL, \cP)$ such that~$\bigcup \cP$ is immersed in~$G$.
\end{lemma}

\begin{proof}
    For every oriented \unit\ \csep~$s$ in~$S$, let~$v_s$ be the separator of~$s$.
    Put~$L := \set{v_s : s \in S}$.
    The linear order on~$S$ induces the linear order~$\leq_L$ on~$L$ that is defined by letting~$v_s\leq_L v_t$ whenever~$s\leq t$.
    Then~$\cL:=(L,\leq_L)$ is order-isomorphic to~$(S, \leq)$, and in particular to~$\QQ\cap [0,1]$.
    
    Let~$x$ and~$y$ be the minimal and maximal element of~$\cL$, respectively.
    Since~$G$ is infinitely edge-connected, there is an infinite system~$\cP''$ of pairwise edge-disjoint~$x$--$y$ paths in~$G$.
    By~\cref{thm:existenceofgrainlines}, we find an~$x$--$y$ grain line~$(\cL', \cP')$ in~$G$ such that every path in~$\cP'$ is also contained in~$\cP''$.
    By the definition of~$x$ and~$y$, we have that, for every oriented \csep~$s = (A,B) \in S$, the vertices~$x$ and~$y$ are contained in~$A$ and~$B$, respectively. 
    Thus, only finitely many paths in~$\cP''$ avoid~$v_s$.
    Moreover, if~$s$ and~$t$ are two oriented \csep s in~$S$ with~$s \leq t$, then, for all but finitely many~$x$--$y$ paths~$P$ in~$\cP''$, the vertex~$v_s$ comes before~$v_t$ on~$P$.
    Hence, $L \subseteq L'$ and~$\leq_L$ agrees with~$\leq_{L'}$ on~$L$.
    Since~$\cL$ is order-isomorphic to~$\QQ \cap [0,1]$, we may additionally assume, by passing to a subsequence of~$\cP'$, that for every $n\ge 1$ and every two vertices~$u$ and~$v$ in~$L \cap L'_{<n}$ with~$u <_L v$, there is an internal vertex~$w$ of~$u P'_n v$ which is contained in~$L$ and satisfies~$u <_L w <_L v$.
    Let us denote this property by~$(\ast)$.

    For every~$n\in\NN$, let~$P_n$ be the path obtained from~$P_n'\in\cP'$ by replacing each~$L$-path~$Q\subset P_n'$ that has an internal vertex with an edge~$e(Q)$ that joins the ends of~$Q$ (this edge need not be an edge in~$G$).
    By~$(\ast)$, the paths in~$\cP:=(P_n)_{n\in\NN}$ are pairwise edge-disjoint, and~$(\cL,\cP)$ is a widly presented grain line.
    Moreover, $\bigcup\cP$ is immersed in~$\bigcup\cP'$: we take the identity on the vertices and edges, except that we map every edge of the form~$e(Q)$ to the path~$Q$.
\end{proof}

\begin{proof}[Proof of \cref{cor:Qorderedgivesfarey}]
    This follows directly from \cref{lemma:fareygraphinwildgrainline} and \cref{lemma:Qorderedsepgivewildgrainline}.
\end{proof}

\subsection{Decomposing along \csep s} \label{sec:decomposealongcsep}

The following key lemma is the final ingredient for the proof of \cref{mainthm:halvedFareygraphistypical}:

\begin{keylemma} \label{lemma:halvedFarey2aleph0decomp}
    Let~$G$ be an infinitely edge-connected graph.
    Then at least one of the following is true:
    \begin{enumerate}
        \item \label{prop:halvedfareygraphinG} The halved Farey graph is immersed in~$G$;
        \item \label{item:2ccon} $G$ is~$2$-\ccon ;
        \item \label{item:csep2ccon} there is a \unit\ \csep~$\{A,B\}$ of~$G$ such that~$G[B]$ is~$2$-\ccon .
    \end{enumerate}
\end{keylemma}

\begin{proof}
    Let us assume that \cref{prop:halvedfareygraphinG} does not hold; in particular, $\Kaleph$ is not immersed in~$G$.
    Thus, by \cref{lemma:nestedsetseparatingeverything}, there is a nested set~$N$ of \unit\ \csep s of~$G$ which is faithful to~$G$.
    
    Let~$\vN$ be the set consisting of all orientations of the \csep s in~$N$, equipped with the usual order~${\le}$ of oriented separations.
    Note that we may assume that~$\vN$ is non-empty: otherwise, $G$ is~$2$-\ccon , and thus \cref{item:2ccon} holds.
    
    We claim that if $(\vN, \leq)$ has no maximal or minimal element, then \cref{item:csep2ccon} holds.
    So suppose that $(A,B)$ is maximal in~$(\vN,\le)$, say (the other case is analogous). Then $G[B]$ is infinitely edge-connected by
    \cref{lemma:splittinglemma}.
    Now if we assume for a contradiction that $G[B]$ is not $2$-\ccon , then there is a \unit\ \csep~$\set{C,D}$ of~$G[B]$ that separates some two vertices~$u$ and~$v$ of~$G[B]$.
    After renaming the sides of $\{C,D\}$, we may assume that $C$ contains the separator~$w$ of~$(A,B)$.
    Then $\set{A \cup C, D}$ is a \unit\ \csep\ of~$G$ which still separates~$u$ and~$v$.
    However, no \csep\ in~$N$ separates $u$ and~$v$, since $u$ and $v$ lie in $B$ and~$(A,B)$ is maximal in~$(\vN,\le)$.
    Hence $N$ is not faithful to~$G$, a contradiction.

    Therefore, we may assume that $(\vN,\le)$ has no maximal or minimal element.
    Hence, using Zorn's lemma, we find an infinite subset~$S'$ of~$\vN$ such that:
    \begin{itemize}
        \item $\leq$ orders~$S'$ linearly;
        \item $S'$ is unbounded in~$\vN$;
        \item $S'$ is an interval of~$\vN$: for every two~$r, t \in S'$ and~$s \in \vN$ with~$r \leq s \leq t$, we have~$s \in S'$.
    \end{itemize}
    Let~$w$ be any vertex of~$G$.
    Since~$N$ is faithful to~$G$, it either does not contain a \csep\ with separator~$w$, or, if~$w$ is a \ccv\ of~$G$, then there is an orientation of the set of all \csep s in~$N$ with separator~$w$ which is a star of separations. So it follows from~$(S', \leq)$ being a linear order that there are at most two oriented \csep s in~$S'$ with separator~$w$.
    Moreover, if there are two such separations~$s_1 < s_2$, then the separator of any oriented \csep~$t$  with~$s_1 < t < s_2$ is a subset of~$\set{w}$, thus equal to~$\set{w}$, so~$s_2$ is the successor of~$s_1$ in~$(S',\le)$.
    Let~$S$ be a subset of~$S'$ obtained by deleting exactly one element of every pair of oriented \csep s in~$S'$ with the same separator.
    Thus, the separators of the oriented \csep s in~$S$ are pairwise distinct.
    
    Let~$R$ be the set of all elements of~$S$ that have a successor in~$(S,\le)$. We claim that~$R$ is infinite.
    To see this, let us assume for a contradiction that~$R$ is finite.
    Then~$(S\setminus R,\le)$ is an infinite set equipped with a dense linear order.
    So we can greedily find a subset of~$S\setminus R$ which is ordered by~$\leq$ like~$\QQ$, and apply \cref{cor:Qorderedgivesfarey} to obtain a halved Farey graph immersed in~$G$, contradicting our assumption that \cref{prop:halvedfareygraphinG} fails.
    Thus, $R$ is infinite.

    By passing to a subset of~$R$, we may assume without loss of generality that~$(R, \leq)$ is order-isomorphic to~$\omega$ or to the inverse of~$\omega$.
    In the following, we will assume that~$(R, \leq)$ is order-isomorphic to~$\omega$, as the other case is analogous.
    Let~$s_0 < s_2 < s_4 < \dots$ be an enumeration of~$R$.
    For every~$n \in \NN$, let~$s_{2n+1}$ be the successor of~$s_{2n}$ in~$S$.
    Moreover, for every~$n\in\NN$ we let~$(A_n, B_n):=s_n$, denote the separator of~$s_n$ by~$v_n$, and we let~$V_0 := A_0$ as well as~$V_{n+1} := B_n \cap A_{n+1}$.
    
    We claim that, for every natural number~$n \geq 1$,
    \begin{enumerate}[label=(\arabic*)]
        \item \label{prop:infconbetweensepvertices} either~$v_{n-1}$ and~$v_{n}$ are equal or~$G[V_{n}]$ is infinitely edge-connected, and
        \item \label{prop:2cconetweensepvertices} if~$n$ is odd then~$v_{n-1}$ and~$v_{n}$ are distinct and~$2$-\ccon\ in~$G[V_{n}]$.
    \end{enumerate}
    Let~$n \geq 1$ be any natural number.
    If~$v_{n-1}$ and~$v_{n}$ are distinct, then applying \cref{lemma:splittinglemma} twice yields that~$G[V_{n}]$ is infinitely edge-connected.
    Thus, \cref{prop:infconbetweensepvertices} holds.
    Now we additionally assume that~$n$ is odd.
    Then~$s_{n}$ is the successor of~$s_{n-1}$ in~$S$; in particular, $s_{n-1}$ and~$s_n$ are distinct. 
    Thus, since~the separators of the oriented \csep s in~$S$ are pairwise distinct, $v_{n-1}$ and~$v_{n}$ are distinct.
    We claim that there is no \unit\ \csep\ of~$G[V_{n}]$ that separates~$v_{n-1}$ and~$v_{n}$.
    Indeed, suppose for a contradiction that there is a \unit\ \csep~$\{A,B\}$ of~$G[V_{n}]$ with separator~$w$ such that~$v_{n-1}\in A\setminus B$ and~$v_n\in B\setminus A$.
    Then~$\{A_{n-1} \cup A, B \cup B_{n}\}$ is a \unit\ \csep\ of~$G$ which separates~$v_{n-1}$ and~$v_n$.
    Since~$N$ is faithful to~$G$, there is a \csep~$s$ with the same separator~$w$ in~$N$ which separates~$v_{n-1}$ from~$v_n$, and
    by \cref{lemma:distinctseparatorsarenested} we have that~$s_{n-1} < s < s_{n}$ after suitably orienting~$s$.
    As~$S'$ is an interval of~$\vN$, we have~$s \in S'$.
    By the construction of~$S$, there is a unique oriented \csep~$t\in S$ with separator~$w$. Since we have~$s_{n-1}<s<s_n$ and since~$w$ is distinct from the separators~$v_{n-1}$ and~$v_n$ of~$s_{n-1}$ and~$s_n$, \cref{lemma:distinctseparatorsarenested} implies~$s_{n-1}<t<s_n$.
    This contradicts the fact that~$s_{n}$ is the successor of~$s_{n-1}$ in~$S$. So \cref{prop:2cconetweensepvertices} holds.
    
    To conclude the proof, we construct an immersion of~$\Kaleph$ in~$G$, as follows.
    For every even natural number~$n\ge 1$, if~$v_{n-1}$ is distinct from~$v_n$ then we pick an infinite system~$\cP_n$ of infinitely many pairwise edge-disjoint~$v_{n-1}$--$v_n$ paths in~$G[V_n]$, which exists by~\cref{prop:infconbetweensepvertices}.
    For every odd natural number~$n\ge 1$, we pick an arbitrary vertex~$u_n$ in~$V_{n} \setminus \{v_{n-1}, v_{n}\}$ and---using that~$G[V_n]$ is~$2$-\ccon\ by~\cref{prop:2cconetweensepvertices}---we greedily pick an infinite system~$\cP_n$ of pairwise edge-disjoint~$v_{n-1}$--$v_{n}$ paths in~$G[V_{n}]$ such that infinitely many of these contain~$u_{n}$ and infinitely many do not. Note that a vertex~$u_n$ can only lie on paths in~$\cP_n$ and not on paths in any other~$\cP_m$.
    Let~$U := \set{u_{n} : n \in 2\NN+1}$.
    Finally, we greedily connect every two distinct vertices~$u_i, u_j \in U$ with a path avoiding all other~$u_n$ and the already used finitely many edges by suitably concatenating paths from the sets~$\cP_m$ with~$m\ge 1$.
    This yields an immersion of~$\Kaleph$ in~$G$ with~$U$ as its set of branch vertices, and thus \cref{prop:halvedfareygraphinG} holds.
\end{proof}

\subsection{The halved Farey graph is everywhere}\label{sec:proofofthm1}

Finally, we are ready to prove \cref{mainthm:halvedFareygraphistypical}:

\begin{proof}[Proof of \cref{mainthm:halvedFareygraphistypical}]
    Let~$G$ be any infinitely edge-connected graph, and let us assume for a contradiction that the halved Farey graph is not immersed in~$G$.
    We will recursively construct infinitely edge-connected subgraphs~$G_0, G_1, \dots$ and~$H_0, H_1, \dots$ of~$G$ such that, for every~$n \in \NN$,
    \begin{enumerate}
        \item\label{PWED} $\set{G_0, \dots, G_n} \cup \set{H_n}$ is a set of pairwise edge-disjoint subgraphs of~$G$,
        \item\label{GHCSEP} if $n\ge 1$ then $\set{V(G_n), V(H_n)}$ is a \csep\ of~$H_{n-1}$ with~$V(G_n)\cap V(H_n)\neq\emptyset$, and
        \item\label{IG} the intersection graph of~$\set{V(G_m) : m \leq n}$ is connected.
    \end{enumerate}
    
    Later, we will continue working with the graphs~$G_n$, while the graphs~$H_n$ only serve as `reservoirs' during the construction of these two sequences, which we explain next.
    
    To get started, suppose that~$n=0$, so we have to define~$G_0$ and~$H_0$.
    By \cref{fact:findcsep}, there is a \csep~$\{A_0,B_0\}$ of~$G$ that separates two vertices of~$G$ minimally.
    By \cref{lemma:splittinglemma}, $G[A_0]$ and~$G[B_0]$ are infinitely edge-connected.
    Hence, $G_0 := G[A_0]$ and~$H_0 := G[B_0] - E_G(A_0 \cap B_0)$ are the desired subgraphs, where~$E_G(A_0 \cap B_0)$ denotes the set of all the finitely many edges of~$G$ that run between the vertices in~$A_0 \cap B_0$.

    Suppose next that~$n\ge 1$ and that we have already chosen~$G_0, \dots, G_{n-1}$ and~$H_0, \dots, H_{n-1}$.
    By \cref{fact:findcsep}, there is a \csep~$\{A_n,B_n\}$ of~$H_{n-1}$ that separates two vertices of~$H_{n-1}$ minimally.
    Since the intersection~$V(G_{n-1})\cap V(H_{n-1})$ is nonempty by \cref{GHCSEP} for~$n-1$, we can rename the sides of~$\{A_n,B_n\}$ so that~$A_n$ meets~$V(G_{n-1})$.
    By \cref{lemma:splittinglemma}, both~$H_{n-1}[A_n]$ and~$H_{n-1}[B_n]$ are infinitely edge-connected.
    Hence~$G_n := H_{n-1}[A_n]$ and~$H_n := H_{n-1}[B_n] - E_{H_{n-1}}(A_n \cap B_n)$ also are infinitely edge-connected.
    \cref{PWED} holds by construction.
    Since~$H_{n-1}$ is infinitely edge-connected, the intersection~$A_n\cap B_n$ is nonempty, so \cref{GHCSEP} holds.
    Finally, \cref{IG} holds because we named the sides of~$\{A_n,B_n\}$ so that~$A_n$ meets~$V(G_{n-1})$.

    In the next step, we will obtain from the sequence~$(G_n)_{n\in\NN}$ of infinitely edge-connected graphs a countably infinite collection~$\cU$
    of pairwise edge-disjoint infinitely edge-connected subgraphs of~$G$ such that
    \begin{enumerate}[label=(\arabic*)]
        \item\label{M2} infinitely many graphs in~$\cU$ are~$2$-\ccon ,
        \item\label{M3} each graph in~$\cU$ has a vertex that is contained in no other member of~$\cU$, and
        \item\label{M4} the intersection graph of~$\cU$ is connected.
    \end{enumerate}

    We construct~$\cU$ from the empty set by adding one or two graphs for each~$n\in\NN$, as follows.
    Given~$n$, we apply \cref{lemma:halvedFarey2aleph0decomp} to~$G_n$, which tells us that~$G_n$ is~$2$-\ccon\ or that there is a unitary \csep~$\{X,Y\}$ of~$G_n$ such that~$G_n[Y]$ is~$2$-\ccon\ (as the halved Farey graph is not immersed in~$G_n\subset G$ by assumption).
    In the former case, we add~$G_n$ to~$\cU$.
    In the latter case, we add~$G_n[X] - E_{G_n}(X \cap Y)$ and~$G_n[Y]$ to~$\cU$.
    
    The~$2$-\ccon\ graphs in~$\cU$ are infinitely edge-connected by definition, and the graphs of the form~$G_n[X]- E_{G_n}(X\cap Y)$ are infinitely edge-connected by \cref{lemma:splittinglemma}.
    While properties \cref{M2} and \cref{M4} follow immediately from the above construction and~\cref{IG},
    property~\cref{M3} follows from~\cref{GHCSEP}:
    For every~$n \in \NN$, \cref{GHCSEP} implies that only finitely many vertices of~$G_n$ occur in any later~$G_m$ with~$m>n$. Hence every~$G_n$ has only finitely many vertices that occur in any~$G_m$ with~$m\neq n$.
    Thus, every graph in~$\cU$ shares only finitely many vertices with other members of~$\cU$.
    As the graphs in~$\cU$ are infinite, \cref{M3} follows.
    
    To obtain a contradiction, we construct an immersion of~$\Kaleph$ in~$G$, as follows.
    In each~$2$-\ccon\ graph in~$\cU$ we select one vertex that lies in no other member of~$\cU$, which is possible by~\cref{M3}. In total, $\aleph_0$ vertices are select in this way, by \cref{M2}.
    The selected vertices will be the branch vertices of the immersion.
    Now \cref{M2}~and~\cref{M4} combined with the infinite edge-connectivity of the members of~$\cU$ allow us to greedily connect every two branch vertices with a path that avoids all other branch vertices and the finitely many edges used by the previously chosen paths.
\end{proof}

\begin{ack}
    We thank the reviewers for valuable comments that greatly improved the presentation of some arguments and fixed a few mistakes.
\end{ack}

\bibliographystyle{amsplain}
\bibliography{references.bib}
\end{document}